\documentclass[12pt]{amsart}

\usepackage{amssymb,latexsym}
\usepackage[enableskew]{youngtab}
\usepackage{verbatim,euscript}
\usepackage{fullpage}
\usepackage[all]{xy}

\usepackage[colorlinks=true,linkcolor=blue,urlcolor=violet,citecolor=magenta]{hyperref}

\input {cyracc.def}
\numberwithin{equation}{section}

\newtheorem{thm}{Theorem}[section]
\newtheorem{lm}[thm]{Lemma}
\newtheorem{cl}[thm]{Corollary}
\newtheorem{prop}[thm]{Proposition}

\theoremstyle{remark}

\theoremstyle{definition}
\newtheorem{ex}[thm]{Example}
\newtheorem{df}{Definition}
\newtheorem*{rema}{Remark}

\newenvironment{proof*}
{\noindent {\sl Proof.}\quad }{\hfill $\square$}

\newcommand {\ah}{{\mathfrak a}}
\newcommand {\be}{{\mathfrak b}}
\newcommand {\ce}{{\mathfrak c}}

\newcommand {\g}{{\mathfrak g}}

\newcommand {\fH}{{\eus H}}

\newcommand {\el}{{\mathfrak l}}

\newcommand {\p}{{\mathfrak p}}

\newcommand {\te}{{\mathfrak t}}
\newcommand {\ut}{{\mathfrak u}}

\newcommand {\slno}{{\mathfrak{sl}}_{n+1}}

\newcommand {\spn}{{\mathfrak{sp}}_{2n}}
\newcommand {\sono}{{\mathfrak{so}}_{2n+1}}
\newcommand {\sone}{{\mathfrak{so}}_{2n}}

\newcommand {\gD}{{\eus D}}

\newcommand {\gH}{{\eus H}}

\newcommand {\esi}{\varepsilon}
\newcommand {\ap}{\alpha}

\newcommand {\lb}{\lambda}

\newcommand {\HW}{\widehat W}
\newcommand {\HV}{\widehat V}
\newcommand {\HP}{\widehat\Pi}
\newcommand {\HD}{\widehat\Delta}

\newcommand {\EE}{{\mathcal E}}

\newcommand {\ck}{{\mathcal K}}
\newcommand {\ckl}{\check{\mathcal K}}
\newcommand {\cku}{\hat{\mathcal K}}
\newcommand {\LL}{{\mathcal L}}

\newcommand {\PP}{{\mathcal P}}

\newcommand {\rk}{{\mathsf{rk\,}}}

\newcommand {\GR}[2]{{\textrm{{\bf #1}}}_{#2}}

\newcommand {\ov}{\overline}
\newcommand {\un}{\underline}

\newcommand {\Ab}{\mathfrak{Ab}}
\newcommand {\Abo}{\overset{o}{\mathfrak{Ab}}}
\newcommand {\AD}{\mathfrak{Ad}}

\newcommand {\beq}{\begin{equation}}
\newcommand {\eeq}{\end{equation}}

\newcommand{\curle}{\preccurlyeq}
\renewcommand{\le}{\leqslant}
\renewcommand{\ge}{\geqslant}

\newcommand{\adn}{{\sf ad}-nilpotent }

\newcommand{\eus}{\EuScript}
\newcommand{\zab}{\mbox{$\esi_{12}$}}
\newcommand{\zac}{\mbox{$\esi_{13}$}}
\newcommand{\zad}{\mbox{$\esi_{14}$}}
\newcommand{\zae}{\mbox{$\esi_{15}$}}
\newcommand{\zaf}{\mbox{$\esi_{16}$}}
\newcommand{\zbf}{\mbox{$\esi_{26}$}}
\newcommand{\zcf}{\mbox{$\esi_{36}$}}
\newcommand{\zdf}{\mbox{$\esi_{46}$}}
\newcommand{\zef}{\mbox{$\esi_{56}$}}

\newcommand{\za}{\mbox{$\esi_1$}}
\newcommand{\zabn}{\mbox{$\ov{\esi}_{12}$}}
\newcommand{\zacn}{\mbox{$\ov{\esi}_{13}$}}
\newcommand{\zadn}{\mbox{$\ov{\esi}_{14}$}}
\newcommand{\zaen}{\mbox{$\ov{\esi}_{15}$}}
\newcommand{\zafn}{\mbox{$\ov{\esi}_{16}$}}

\newcommand{\zbcn}{\mbox{$\ov{\esi}_{23}$}}
\newcommand{\zbdn}{\mbox{$\ov{\esi}_{24}$}}
\newcommand{\zben}{\mbox{$\ov{\esi}_{25}$}}
\newcommand{\zbfn}{\mbox{$\ov{\esi}_{26}$}}

\newcommand{\zcdn}{\mbox{$\ov{\esi}_{34}$}}
\newcommand{\zcen}{\mbox{$\ov{\esi}_{35}$}}
\newcommand{\zcfn}{\mbox{$\ov{\esi}_{36}$}}

\newcommand{\zden}{\mbox{$\ov{\esi}_{45}$}}
\newcommand{\zdfn}{\mbox{$\ov{\esi}_{46}$}}

\newcommand{\zefn}{\mbox{$\ov{\esi}_{56}$}}

\newfam\Bbbfam\newfam\eufam\newfam\eusfam%
\font\Bbbfont=msbm10 scaled 1200%
\font\Bbbsmallfont=msbm8%
\textfont\Bbbfam=\Bbbfont\scriptfont\Bbbfam=\Bbbsmallfont%%

\begin{document}
\setlength{\parskip}{2pt plus 4pt minus 0pt}
\hfill {\scriptsize February 18, 2012} 
\vskip1.5ex

\title[Two covering polynomials]{Two covering polynomials of a finite 
poset, with applications to root systems and {\sf ad}-nilpotent ideals}
\author{Dmitri I. Panyushev}
%\thanks{This research was supported in part by 
\address[]{Independent University of Moscow,
Bol'shoi Vlasevskii per. 11, 119002 Moscow, \ Russia
\hfil\break\indent
Institute for Information Transmission Problems, B. Karetnyi per. 19, Moscow 
127994}
\email{panyushev@iitp.ru}
\keywords{Root system, ad-nilpotent ideal, graded poset, Hasse diagram}
\subjclass[2010]{06A07, 17B20, 20F55}
\begin{abstract}
We introduce two polynomials (in $q$) associated with a finite poset $\PP$ that encode some information on the covering relation in $\PP$. If $\PP$ is a distributive lattice, and hence 
$\PP$ is isomorphic to the poset of dual order ideals in a poset $\LL$,
then these polynomials coincide and the coefficient of $q$ equals the number of $k$-element 
antichains in $\LL$. In general, these two covering polynomials are different, and 
we introduce a deviation polynomial of $\PP$, which measures the difference between 
these two. We then compute all these polynomials in the case, where $\PP$ is one of 
the posets associated with an irreducible root system. 
These are 1) the posets of positive roots, 2) the poset of \adn  ideals,  and 
3) the poset of Abelian ideals.
\end{abstract}
\maketitle

\section*{Introduction}

\noindent
In this note, we associate two polynomials  with a 
finite poset $\PP$, study their properties, and determine these polynomials 
for some interesting 
posets related to root systems. Specifically, we mean the poset of positive roots,
poset of {\sf ad}-nilpotent ideals, and poset of abelian ideals, see definitions below.
%considered in \cite{rodstv}. 

Consider two statistics $\kappa$ and $\iota$ on $\PP$. 
By definition, $\kappa(x)$, $x\in\PP$, is the number of elements of $\PP$
that are covered by $x$,
and $\iota(x)$ is the number of elements that cover $x$.
The generating function associated with $\kappa$ (resp. $\iota$) is called the 
{\it upper\/} (resp. {\it lower}) {\it covering polynomial\/} of $\PP$. That is,
$%%\displaystyle 
\cku_{\PP}(q)=\sum_{x\in\PP} q^{\kappa(x)}$
and 
$%%\displaystyle 
\ckl_{\PP}(q)=\sum_{x\in\PP} q^{\iota(x)}$.
The upper covering polynomial, $\cku_\PP$, has briefly been considered, 
without the adjective `upper', in
\cite[Section\,5]{rodstv}. 
It immediately follows from the definition that 
$\cku_{\PP}(1)=\ckl_{\PP}(1)$ and
$\cku'_{\PP}(1)=\ckl'_{\PP}(1)$. The last equality stems from the
observation that both values equal the number of edges in
the Hasse diagram of $\PP$. Consequently, 
$\cku_\PP(q)-\ckl_\PP(q)=(q-1)^2 \gD_\PP(q)$ for some polynomial 
$\gD_\PP$, which is called the {\it deviation polynomial\/} of $\PP$.

In Section~\ref{odin}, we begin with  some simple observations on these polynomials and then prove that $\cku_\PP\equiv\ckl_\PP$ if $\PP$ is 
a distributive lattice (see Theorem~\ref{ravno1}). For a distributive lattice $\PP$, the common covering polynomial is denoted by $\mathcal K_\PP$.
% or admits an order-reversing involution.

Let $\Delta$ be an irreducible root system, $\Delta^+$ a subset of
positive roots, and $\Pi\subset\Delta^+$  the set of simple roots.
If $\Delta$ is reduced, then $\g$ is the corresponding
simple Lie algebra, with fixed Borel subalgebra $\be$ corresponding to $\Delta^+$.
We determine $\cku_\PP$, $\ckl_\PP$, and $\gD_\PP$ in the following  cases:

1)  $\PP=\Delta^+$, equipped with the standard root order, see Section~\ref{dva};

2)  $\PP$ is $J^*(\Delta^+)$ or $\PP=J^*(\Delta^+\setminus\Pi)$, where 
$J^*(\LL)$ stands for the poset of dual order ideals in $\LL$. 
%the poset of dual order ideals in $\Delta^+ \setminus\Pi$. 
In the Lie algebra case, $J^*(\Delta^+)$ is isomorphic to the poset of 
{\sf ad}-{\it nilpotent ideals\/} 
of $\be$, i.e., $\be$-ideals in $\ut=[\be,\be]$, the nilradical of $\be$;
and $J^*(\Delta^+\setminus\Pi)$ is isomorphic to the poset
of $\be$-ideals in $[\ut,\ut]$.
These two posets are also denoted by $\AD$ and $\AD_0$, respectively; 

3) $\PP=\Ab$, the subposet of $\AD$ that consists of the abelian ideals of $\be$.
(An {\it abelian ideal\/} of $\be$ is a subspace $\ce\subset\be$ such that $[\be,\ce]\subset \ce$
and $[\ce,\ce]=0$.)
\\[.8ex]
Whenever we wish to stress that these (po)sets depend on $\g$, we write $\Delta(\g)$,
$\AD(\g)$, etc.

Let us briefly describe our results. 
For $\Delta^+$, we prove that $\deg\cku_{\Delta^+}=\deg\ckl_{\Delta^+}\le 3$
and the coefficients of $q^3$ in $\cku_{\Delta^+}$ and $\ckl_{\Delta^+}$ are equal.
This implies that the deviation polynomial is a constant; namely,
$\gD_{\Delta^+}(q)\equiv \rk\Delta-1$. This also includes the only non-reduced irreducible
root system $\GR{BC}{n}$. We show  that the posets $\Delta^+(\GR{BC}{n})$,  
$\Delta^+(\GR{B}{n+1})\setminus \Pi$, and  $\Delta^+(\GR{C}{n+1})\setminus \Pi$
are isomorphic, which allows  to reduce many problems on $\GR{BC}{n}$ to $\GR{B}{n+1}$ or  $\GR{C}{n+1}$.

Since  $\AD$ and $\AD_0$ are distributive lattice, we have only two different covering polynomials for them. The polynomial $\ck_\AD$
%%$\cku_{\AD}\equiv\ckl_{\AD}$  and \cku_{\AD_0}\equiv\ckl_{\AD_}$.
appeared earlier under various guises in different theories in
\cite{ath1,bessis,duality,charney}. The coefficients of $\ck_\AD$ are the {\it generalised Narayana numbers}.
{A posteriori}, it is known  that $\cku_{\AD}$ is palindromic, but no general explanation
is available in the context of \adn ideals.
Explicit formulae for $\ck_{\AD_0}$ show that this polynomial is not always palindromic.
(Our computation of  $\ck_{\AD_0(\GR{D}{n})}$ 
%for series $\GR{D}{n}$ %and $\GR{E}{n}$ 
relies  on the conjectural relationship between the coefficients of
$\ck_{\AD_0}(q)$ and the $\eus F$-triangle introduced by F.\,Chapoton \cite{chap04}, see
Section~\ref{tri}.)
However, the ratio $\cku'(1)/\cku(1)$ is determined by similar rules in both cases.
We notice that 
\[
\frac{\ck'_\AD(1)}{\ck_\AD(1)}=%\frac{n}{2}=
\frac{\#(\Delta^+)}{h}\quad {\text {and}} \quad 
\frac{\ck'_{\AD_0}(1)}{\ck_{\AD_0}(1)}=
%\frac{n}{2}\cdot \frac{h-2}{h-1}=
\frac{\#(\Delta^+\setminus\Pi)}{h-1}\ ,
\]
where $h$ is the Coxeter number of $\Delta$. %and $n$ is the rank of $\Delta$.
The first equality stems from the fact that $\ck_\AD$ is palindromic, 
of degree $n=\rk\Delta $ (although the fact that $\ck_\AD$ is palindromic is not explained yet). 
The reason for the validity of the second one is totally unclear.

The most interesting case is that of abelian ideals. Here the upper and lower
covering polynomials are usually different. The reason is that although 
$\Ab$ is a meet semilattice, it is a distributive lattice if and only if 
$\Delta$ is of type $\GR{C}{n}$ or $\GR{G}{2}$.
We develop some general theory for computing covering polynomials,
which is based on a bijection between the abelian ideals and the minuscule
elements of the affine Weyl group of $\Delta$. 
Let $I\subset\Delta^+$ be an abelian ideal. Using the minuscule element
corresponding to  $I$, we define the {\it shift
vector\/} $\mathbf k_I=(k_0,k_1,\ldots,k_n)$, where $k_i\in\{-1,0,1,2\}$, and
prove that $\kappa(I)=\#\{ i\mid k_i=-1\}$ and $\iota(I)=\#\{j\mid k_j=1\}$.
We then describe a recursive procedure for computing all $\mathbf k_I$ starting
from $I=\varnothing$. The procedure basically asserts that if $k_i=1$, then
$\mathbf k_I$ can be replaced with $\mathbf k_I - \{ \text{the $i$-th column of the extended Cartan matrix of $\Delta$} \}$,
see Section~\ref{genera} for details.

We also present a method of calculation of $\ckl_\Ab$, which exploits
the canonical mapping of $\Ab\setminus\{\varnothing\}$ onto the
set of long positive roots \cite{imrn}. For explicit computations with exceptional root
systems, we use the general equalities
$\cku_\Ab(1)=2^n$ \cite{cp1,ko1} and $\cku'_\Ab(1)=(n+1)2^{n-2}$ \cite{rodstv};
our calculations in the classical cases exploit standard matrix presentations
of these Lie algebras and counting certain Ferrers diagrams.

Our computations show that, for many natural posets, the coefficients of $\gD_\PP$
are of the same sign. This includes $\Ab$, $\Delta^+$, $\Delta^+\cup\{0\}$,
$\Delta^+\setminus\Pi$.
%%Our results show that the coefficients of $\gD_{\Ab}$ are non-positive. 
It is likely that there could exist some general
conditions on $\PP$ guaranteeing that $\gD_\PP$ has the coefficients of the
same sign. In Section~\ref{some_spec}, we propose a condition of such sort.

{\small
{\bf Acknowledgements.} This work was done during my visits to
the Max-Planck-Institut f\"ur Mathematik (Bonn), and I thank the Institute for
the hospitality and inspiring environment. I wish to thank F.\,Chapoton for sending
me his unpublished notes and R.\,Stanley for drawing my attention to work of Dilworth \cite{dilw}.
}

%%%%%%%%%%%%%%%%   Section 1

\section{Definition and basic properties}
\label{odin}

\noindent
Let $(\PP, \curle)$ be a finite poset. 
Write $\fH(\PP)$ for the Hasse diagram of $\PP$ and
$\EE(\PP)$ for the set of edges
of $\fH(\PP)$. We regard $\fH(\PP)$ as a digraph; if $x$ covers $y$
($x,y\in\PP$), then the edge $(x,y)$ is depicted as $y\to x$ and we say that  $(x,y)$ {\it originates\/} in $y$ and {\it terminates\/} in $x$.

For any $x\in\PP$, let $\kappa(x)$
be the number of $y\in \PP$ such that $y$ is covered by $x$,
and $\iota(x)$  the number of $y\in \PP$ such that $y$  covers $x$.
In other words, $\kappa(x)$ (resp. $\iota(x)$) 
is the number of edges in $\EE(\PP)$ that terminates (resp. originates)  in $x$.
We define two polynomials that encode some
properties of the covering relation in $\PP$. 

\begin{df}  \label{cover_pol} \leavevmode\par
{\sf (i)} \ The {\it upper covering polynomial\/} of $\PP$ is
$\displaystyle \cku_{\PP}(q)=\sum_{x\in\PP} q^{\kappa(x)}$;

{\sf (ii)} \ The {\it lower covering polynomial\/} of $\PP$ is
$\displaystyle \ckl_{\PP}(q)=\sum_{x\in\PP} q^{\iota(x)}$;
\end{df}%

\noindent
It follows that $\cku_\PP(0)$ (resp. $\ckl_\PP(0)$)
is the number of the minimal (resp. maximal) elements of $\PP$. 
In general, these polynomials are different; they may even have different 
degree. However, one readily deduces from the definition that 
\[
\cku_{\PP}(q)\vert_{q=1}=\ckl_{\PP}(q)\vert_{q=1}=
\#\PP \quad  \text{and} \quad
\frac{d}{dq}\cku_\PP(q)\vert_{q=1}=\frac{d}{dq}\ckl_\PP(q)\vert_{q=1}=
\#\EE(\PP) \ .
\]
Hence $\cku_\PP(q)-\ckl_\PP(q)=(q-1)^2 \gD_\PP(q)$ for some polynomial
$\gD_\PP$. We will say that $\gD_\PP$ is the {\it deviation polynomial\/} 
of $\PP$.   The following is straightforward. 
\begin{lm}   \label{prod} 
We have the following properties%\leavevmode\par
\begin{itemize}
\item[\sf (i)] \  If\, $\PP=\PP_1+\PP_2$, then 
$\cku_\PP=\cku_{\PP_1}+\cku_{\PP_2}$, and likewise for $\ckl$ and $\gD$;
\item[\sf (ii)] \ If\, $\PP=\PP_1\times\PP_2$, then 
$\cku_\PP=\cku_{\PP_1}\cku_{\PP_2}$, $\ckl_\PP=\ckl_{\PP_1}\ckl_{\PP_2}$,
and $\gD_\PP=\cku_{\PP_1}\gD_{\PP_2}+\ckl_{\PP_2}\gD_{\PP_1}=
\ckl_{\PP_1}\gD_{\PP_2}+\cku_{\PP_2}\gD_{\PP_1}$.
\end{itemize}
\end{lm}

\noindent We are going to investigate how properties of $\PP$ are reflected 
in $\cku_\PP, \ckl_\PP, \gD_\PP$.

\begin{thm}   \label{ravno1}
Let\/ $\PP$ be a distributive lattice. Then $\cku_\PP= \ckl_\PP$.
More precisely, if $\PP\simeq J(\LL)$, then the coefficient
of $q^k$ equals the number of $k$-element antichains in $\LL$.
\end{thm}
\begin{proof}
By Birkhoff's theorem for finite distributive lattices, $\PP$ is 
isomorphic to the poset of order ideals of a unique poset $\LL$,
i.e., $\PP\simeq J(\LL)$, see e.g. 
\cite[Theorem\,3.4.1]{rstan1}. 
If $I$ is an order ideal of $\LL$, then the set of maximal elements of
$I$, $\max(I)$, is an antichain of $\LL$. And the same is true for the set of minimal 
elements of $\LL\setminus I$, $\min(\LL\setminus I)$.
It easily follows from Definiton~\ref{cover_pol} that
regarding $I$ as an element of $\PP$ we have
$\iota(I)=\# \max(I)$ and $\kappa(I)=\#\min(\LL\setminus I)$.
Conversely, each antichain in $\LL$ occurs as both $\max(I)$ and 
$\min(\LL\setminus J)$ for suitable order ideals $I,J$.
This means that both covering polynomials essentially count all
the antichains of $\LL$ with respect to their cardinality.
\end{proof}

{\it Remark.} 
More generally, the equality $\cku_\PP= \ckl_\PP$ holds if $\PP$ is a modular lattice \cite{dilw}. This result of Dilworth  is also discussed  
in \cite[Ex.\,3.38.5]{rstan1}. 
%This fact was communicated to me by R.\,Stanley.

Below, it will be more convenient for us to think of a 
distributive lattice as the poset of {\it dual order\/} (= {\it upper}) ideals.
The distributive lattice of upper ideals of a poset $\LL$ is denoted by $J^*(\LL)$.
Then $\LL$ is being restored as the set of meet-irreducibles in $J^*(\LL)$.
If $I\in J^*(\LL)$, then $\kappa(I)=\# \min(I)$ and
$\iota(I)=\#\max (\LL\setminus I)$.  Note also that  the posets $J^*(\LL)$ and $J(\LL^{op})$ 
are canonically isomorphic.

The following is clear:

\begin{prop}   \label{ravno2}
If\/ $\PP$ admits an order-reversing bijection, then 
$\cku_\PP=\ckl_\PP$.
\end{prop}
%\begin{proof*} If $\omega: \PP\to \PP$ is an order-reversing bijection, then
%$\kappa(x)=\iota(\omega(x))$ for any $x\in\PP$. \end{proof*}

\begin{ex} 1$^o$. Let $(W,S)$ be a finite Coxeter group.
Consider $W$ as poset under the Bruhat-Chevalley ordering `$\le$'.
It is easily seen that $W$ is not a lattice. But the mapping 
$w\mapsto ww_0$, where $w_0\in W$ is the longest element, yields an
order-reversing bijection of $(W,\le)$. Hence $\cku_{W}(q)=\ckl_W(q)$.
More generally, such an equality also holds for the poset $W/W_J$, where $J\subset S$
and $W_J$ is the corresponding parabolic subgroup of $W$.

2$^o$. Let $\PP$ be an arbitrary poset and $\PP^{op}$ the opposite poset.
Then $\cku_\PP=\ckl_{\PP^{op}}$ and $\ckl_\PP=\cku_{\PP^{op}}$. Hence
$\gD_{\PP}=-\gD_{\PP^{op}}$. It then follows from Lemma~\ref{prod}
that $\gD_{\PP\times\PP^{op}}=0$. One may also notice that
$\PP\times\PP^{op}$ admits an order-reversing involution.
\end{ex}%

\begin{ex}    \label{ex:transform}
As we have shown, if $\PP=J^*(\mathcal L)$, then $\gD_\PP= 0$.
%Let us look at the effect of two simple transformations of a distributive lattice.
Let $\hat 0$ and $\hat 1$ denote the maximal and minimal element of
$\PP$, respectively. Set $\PP'=\PP\setminus \{\hat 0\}$
and $\PP''=\PP\setminus\{\hat 1\}$. It is easy to realise the
effect of these procedures for both covering polynomials.
%$\cku_\PP$ and $\ckl_\PP$. 
%Let $m$ (resp. $l$) be the number of maximal (resp. minimal) elements, then
%\begin{itemize}

\textbullet\quad 
Let $m$  be the number of maximal  elements  of $\LL$. Then 
$\cku_{\PP'}(q)=\ck_\PP(q)-mq+m-1$ and $\ckl_{\PP'}(q)=\ck_\PP(q)-q^m$.
Hence $\cku_{\PP'}(q)-\ckl_{\PP'}(q)=q^m-mq+(m-1)$ and
$\gD_{\PP'}(q)=q^{m-2}+2q^{m-3}+\ldots +(m-2)q+(m-1)$. In particular, $\gD_{\PP'}\equiv 0$
if and only if $m=1$. This is not surprising, because if $\LL$ has a unique maximal element,
say $a_1$, then $\PP\setminus \{\hat 0\}\simeq J^*(\LL\setminus \{a_1\})$.

\textbullet\quad 
Since  $J^*(\LL)\setminus \{\hat 1\} \simeq J(\LL^{op})\setminus \{\hat 0\}$,
the formulae for $\PP''$ are similar.
Only the roles of $\cku$ and $\ckl$ are reversed, and in place of $m$ we need 
the number, say $l$, of minimal  elements  of $\LL$. 
%Then $\cku_{\PP''}=\ck_\PP-q^l$ and $\ckl_{\PP''}=\ck_\PP-l(q-1)-1$.
Therefore, %$\cku_{\PP''}-\ckl_{\PP''}=-q^l+lq-(l-1)$ and
$\gD_{\PP''}(q)=-\bigl(q^{l-2}+2q^{l-3}+\ldots +(l-2)q+(l-1)\bigr)$.
%In particular, $\gD_{\PP''}\equiv 0$ if and only if $l=1$.
%In both cases, we obtain polynomials having all nonzero coefficients of the same sign.
\end{ex}

%$3^o$. 
From the definition of $\gD_\PP$, it follows that
$\gD_\PP(1)=\frac{1}{2}(\cku''_\PP(1)-\ckl''_\PP(1))$.
Therefore
\[ 
\gD_\PP(1)=\frac{\#\{(x,y_1,y_2)\in \PP^3 \mid y_1{\to} x,\ y_2{\to} x\}-
  \#\{(x_1,x_2,y)\in \PP^3 \mid y{\to} x_1,\ y{\to} x_2\}}{2},
\]
where it is assumed that $y_1\ne y_2$ and $x_1\ne x_2$.
In other words, 
\[
2\gD_\PP(1)=
\# \left\{ \text{\begin{picture}(32,20)(0,7)
\put(16,20){\circle{4}}
\put(1,5){\circle{4}}\put(31,5){\circle{4}}
\put(28,7){\vector(-1,1){11}}
\put(2,7){\vector(1,1){11}}
\end{picture}}
\right\}-
\#\left\{\text{\begin{picture}(32,20)(0,7)
\put(16,5){\circle{4}}
\put(1,20){\circle{4}}\put(31,20){\circle{4}}
\put(14,7){\vector(-1,1){11}}
\put(18,7){\vector(1,1){11}}
\end{picture}}
\right\},
\]
the difference between the number of two types of configurations in $\gH(\PP)$.
These configurations are said to be $\wedge$-triples and $\vee$-triples,
respectively.
Using this interpretation, one obtains the following result.

\begin{prop}
Suppose $\tilde\PP$ is a distributive lattice and $\PP\subset\tilde\PP$
a subposet such that if $I\in\PP$ and $I'\curle I$ ($I'\in\tilde\PP$), then $I'\in\PP$.
Then $\gD_\PP(1)\le 0$. Furthermore, $\gD_\PP(1)=0$ if and only if
$\PP$ is a distributive lattice if and only if $\gD_\PP\equiv 0$.
\end{prop}
\begin{proof}
Here each $\wedge$-triple can be completed to a {\it diamond} inside 
$\PP$, i.e., the configuration of
the form \ `$\Diamond$'. This provides an injection of the set of $\wedge$-triples
to the set of $\vee$-triples. If this is a bijection, i.e., each $\vee$-triple can be
included in a diamond, then $\PP$ has a unique maximal element. Hence $\PP$
is a distributive lattice and $\gD_\PP=0$.
\end{proof}%

\noindent
In the following sections, we consider the polynomials $\cku_\PP$,
$\ckl_\PP$, and $\gD_\PP$ for the posets described in the Introduction.
%associated  with systems of positive roots. 

%%%%%%%%%%%%%%%%   Section 2

\section{Covering polynomials for the root systems}
\label{dva}
\setcounter{equation}{0}

\noindent
Our main reference for root systems and their properties is~\cite{bour}. 
Let $\Delta$ be a root system in an $n$-dimensional real euclidean
vector space $V$. Choose a subsystem of positive roots $\Delta^+$ with 
the corresponding set of simple roots $\Pi=\{\ap_1,\ldots,\ap_n\}$.
Write $\theta$ for the highest root in $\Delta^+$ and $h$ for the
{\it Coxeter number}. The standard root order `$\curle$' in $\Delta^+$ is determined by the condition that
$\gamma$ covers $\mu$ if and only if $\gamma-\mu\in \Pi$.
Then $\theta$ is the unique maximal element of $\Delta^+$.
Our goal is to compute both covering polynomials for $(\Delta^+, \curle)$. 
In view of Lemma~\ref{prod}, it suffices to consider the
irreducible root systems. 

If $\Delta$ is reduced and irreducible, then $\g$ is the corresponding simple Lie algebra.
If $\theta=\sum_{i=1}^n m_i\ap_i$, then in this case one also has
$\sum_i m_i=h-1$ \cite[Ch.\,VI\ \S\,1.11, Prop.\,31]{bour}.

In what follows, $[q^m]\mathcal F$ stands for
the coefficient of $q^m$ in the polynomial ${\mathcal F}(q)$.

\begin{thm}  \label{h-roots}
Let $\Delta$ be an irreducible root system of rank $n$. Then
\begin{itemize} 
\item[\sf (i)] \  $\deg \cku_{\Delta^+}\le 3$ and $\deg \ckl_{\Delta^+}\le 3$;
\item[\sf (ii)] \  $[q^3] \cku_{\Delta^+}=[q^3]\ckl_{\Delta^+}$;
\item[\sf (iii)] \ if $\Delta$ is reduced, then $[q] \cku_{\Delta^+}=h-1-n$;
\item[\sf (iv)] \ 
if $\Delta$ is simply-laced, then $[q] \cku_{\Delta^+}=
[q^3]\cku_{\Delta^+}$;
\item[\sf (v)] \  $\gD_{\Delta^+}(q)\equiv n-1$.
\end{itemize}
\end{thm}
\begin{proof} We provide a uniform proof for parts (i) and (ii)
only in the simply-laced case. The remaining cases (including the
non-reduced root system $\GR{BC}{n}$, see below)
can be handled in a case-by-case fashion.

(i)  For $\cku_{\Delta^+}$, 
one has to show that there are at most 3 simple roots that can be
subtracted from a positive root.
%%First, assume that $\Delta$ is simply laced. 
Suppose $\gamma\in \Delta^+$
and $\gamma-\ap_i\in \Delta^+$ for $\ap_i\in\Pi$ and $i=1,2,\dots,k$.
If $(\ap_1,\ap_2)\ne 0$, then these two adjacent simple roots generate the root 
system of type $\GR{A}{2}$. Furthermore,
$\gamma$ is the highest weight in the adjoint 
$\GR{A}{2}$-module inside $\g$. 
Therefore the weight $\gamma-\ap_1-\ap_2$ has multiplicity
two. This is only possible if $\gamma=\ap_1+\ap_2$ and hence $k=2$. 
This also proves that if $k\ge 3$, then all simple roots that can be subtracted
from $\gamma$ are pairwise orthogonal. In this situation, it was shown in
\cite[Corollary\,3.3]{aif99} that $k\le 3$.

The argument for $\ckl_{\Delta^+}$ is similar.

(ii) \ Suppose that $\kappa(\gamma)=3$, and let $\ap_1,\ap_2,\ap_3$ be the
corresponding simple roots. As is shown in part (i), these roots are 
pairwise orthogonal. Therefore $\gamma-\ap_1-\ap_2-\ap_3\in \Delta^+$ and
the mapping $\gamma\mapsto (\gamma-\ap_1-\ap_2-\ap_3)$ sets up a bijection
between $\{\gamma\in \Delta^+\mid \kappa(\gamma)=3\}$ and
$\{\gamma\in \Delta^+\mid \iota(\gamma)=3\}$.

(iii) \ If $\gamma=\sum_j a_j\ap_j\in \Delta^+$, then $[\gamma:\ap_j]:=a_j$ is called the 
$\ap_j$-{\it height\,} of $\gamma$. Suppose that $\kappa(\gamma)=1$, i.e., there is a unique 
$\ap_i\in \Pi$ such that $\gamma-\ap_i\in \Delta^+$.  Let $\el_i$ denote the semisimple subalgebra of $\g$ whose set of simple roots if $\Pi\setminus \{\ap_i\}$. All the roots with
a fixed $\ap_i$-height form the set of weights of a \un{sim}p\un{le} $\el_i$-submodule inside 
$\g$ \cite[2.1]{jos}. Therefore, for each value $\ge 2$ of $\ap_i$-height, there is a unique root $\gamma$
with such property.  Hence, the number of positive roots $\gamma$ 
such that $\ap_i$ is the only simple root that can be subtracted from $\gamma$ equals
$m_i-1$.  Thus, 
$\#\{\gamma\in\Delta^+\mid \kappa(\gamma)=1\}=\sum_{i=1}^n(m_i-1)=h-1-n$.

(iv) \ As is well known,
the number of positive roots is $\cku_{\Delta^+}(1)=nh/2$.
%where $h$ is the Coxeter number of $\Delta$. 
By  \cite[Theorem\,1.1]{rodstv},
the number of edges of $\gH(\Delta^+)$ equals $n(h-2)$ in the simply laced case.
That is, $\cku'_{\Delta^+}(1)=n(h-2)$.
Writing $\cku_{\Delta^+}(q)=n +aq+bq^2+cq^3$ and using the above two 
equalities, we obtain $a=c$.

(v) \ By parts (i) and (ii), $\deg(\cku_{\Delta^+}-\ckl_{\Delta^+})\le 2$.
It is also clear that $[q^0](\cku_{\Delta^+}-\ckl_{\Delta^+})=n-1$.
Since $(q-1)^2$ divides this polynomial, the quotient must be $n-1$.
\end{proof}

Formulae of the theorem, together with known values $\cku_{\Delta^+}(1)$ and
$\cku'_{\Delta^+}(1)$, allow us to write down a closed formula for $\cku_{\Delta^+}$
in the simply-laced case:

\begin{cl}
In the {\bf A-D-E} case,  we have 
\[ \cku_{\Delta^+}(q)=1+(h-1-n)q+ 
(\frac{nh}{2}+n-2h+2)q^2+(h-1-n)q^3 .
\]
\end{cl}

\noindent
It is not hard to compute covering polynomials for the posets $\Delta^+$
in the remaining cases, see Table~\ref{table:roots}.
%%Here is the answer:

\begin{table}[htb]  
\begin{tabular}{cll}
$\Delta$ & $\cku_{\Delta^+}(q)$ & $\ckl_{\Delta^+}(q)$ \\ \hline
$\GR{A}{n}$ & $n +\genfrac{(}{)}{0pt}{}{n}{2}q^2$ {\rule{0pt}{2.5ex}}
& $1+(2n-2)q+\genfrac{(}{)}{0pt}{}{n-1}{2}q^2$ 
\\
$\GR{B}{n},\GR{C}{n}$ & $n + (n-1)q+ (n-1)^2q^2$ {\rule{0pt}{2.4ex}}
& $1 + (3n-3)q+ (n-1)(n-2)q^2$ 
\\
$\GR{BC}{n}$ & $n + nq+ n(n-1)q^2$ {\rule{0pt}{2.4ex}}
& $1 + (3n-2)q+ (n-1)^2q^2$ 
\\
$\GR{D}{n}$ & $n{+}(n{-}3)q{+}(\genfrac{(}{)}{0pt}{}{n}{2}{+}
\genfrac{(}{)}{0pt}{}{n-3}{2})
q^2 {+}(n{-}3)q^3$ {\rule{0pt}{2.3ex}} 
& $1{+}(3n{-}5)q{+}(\genfrac{(}{)}{0pt}{}{n-1}{2}{+}
\genfrac{(}{)}{0pt}{}{n-3}{2})
q^2{+}(n{-}3)q^3$ 
\\
$\GR{E}{6}$ & $6+\phantom{1} 5q+20q^2+\phantom{1} 5q^3$  {\rule{0pt}{2.4ex}}
& $1+15q+15q^2+\phantom{1} 5q^3$ 
\\
$\GR{E}{7}$ & $7+10q+36q^2+10q^3$ & $1+22q+30q^2+10q^3$ 
\\
$\GR{E}{8}$ & $8+21q+70q^2+21q^3$ & $1+35q+63q^2+21q^3$ 
\\
$\GR{F}{4}$ & $4+\phantom{1} 7q+12q^2+\phantom{21} q^3$ &
$1+13q+\phantom{1} 9q^2+\phantom{21} q^3$ 
\\
$\GR{G}{2}$ & $2+\phantom{1} 3q+\phantom{12} q^2$ & $1+\phantom{1} 5q$
\\ 
\hline
\end{tabular}
\vskip1ex 
\caption{The  covering polynomials for the root systems}
\label{table:roots}
\end{table}%
%\begin{rema} The degree of polynomials $\cku$ and $\ckl$ equals 3 if and only
%if the Dynkin diagram of $\Delta$ has a branching node.
%\end{rema}
\noindent
Let $\esi_1,\ldots,\esi_n$ be an orthonormal basis for $V$.
Recall that 
\[
 \Delta(\GR{B}{n})=\{ \text{$\pm\esi_i\pm\esi_j$ ($1\le i < j\le n$),
 $\pm \esi_i$ ($1\le i\le n$)} \},
\]
\[
 \Delta(\GR{C}{n})=\{ \text{$\pm\esi_i\pm\esi_j$ ($1\le i < j\le n$),
 $\pm 2\esi_i$ ($1\le i\le n$)} \}
\]
and the unique non-reduced irreducible root system
$\GR{BC}{n}$ is the union of these two, i.e.,
\[
 \Delta(\GR{BC}{n})=\{ \text{$\pm\esi_i\pm\esi_j$ ($1\le i < j\le n$),
 $\pm \esi_i$, $\pm 2\esi_i$ ($1\le i\le n$)} \}.
\]
%The simple roots of $\Delta(\GR{BC}{n})$ are $\ap_i=\esi_i-\esi_{i+1}$, $1\le i\le n-1$, 
%and $\ap_n=\esi_n$.
The following observation reduces many questions about $\GR{BC}{n}$
to $\GR{B}{n+1}$ or $\GR{C}{n+1}$.

\begin{lm}  \label{poset-BC}
The poset of positive roots for $\GR{BC}{n}$ is isomorphic
to the subposet of non-simple positive roots for $\GR{B}{n+1}$
or\/ $\GR{C}{n+1}$.
The posets $\Delta^+(\GR{B}{n+1})$ and $\Delta^+(\GR{C}{n+1})$ are isomorphic.
\end{lm}
\begin{proof*}
An order-preserving bijection between $\Delta^+(\GR{BC}{n})$ and either
$\Delta^+(\GR{B}{n+1})\setminus\Pi$ or
$\Delta^+(\GR{C}{n+1})\setminus\Pi$ is given as follows:
\[
\begin{array}{rcccll}
     \GR{C}{n+1} &      &  \GR{BC}{n} &      & \GR{B}{n+1}\\
\esi_i-\esi_{j+1} & \xymatrix{{}&{}\ar@{|->}[l]}
%%\gets 
& \esi_i-\esi_j  &  \xymatrix{{}&{}\ar@{<-|}[l]}
%%\mapsto 
& \esi_i-\esi_{j+1} 
                                                   & (1\le i < j\le n)\\
\esi_i+\esi_j     & \xymatrix{{}&{}\ar@{|->}[l]} & \esi_i+\esi_j  & \xymatrix{{}&{}\ar@{<-|}[l]} & \esi_i+\esi_{j+1} 
                                                   & (1\le i < j\le n)\\
\esi_i+\esi_{n+1} & \xymatrix{{}&{}\ar@{|->}[l]} &  \esi_i        & \xymatrix{{}&{}\ar@{<-|}[l]} &  \esi_i 
                                                   & (1\le i\le n) \\
      2\esi_i     & \xymatrix{{}&{}\ar@{|->}[l]} &  2\esi_i       & \xymatrix{{}&{}\ar@{<-|}[l]} &  \esi_i+\esi_{i+1} 
                                                   & (1\le i\le n) 
\end{array}
\]
It is easily seen that this extends to an isomorphism between 
$\Delta^+(\GR{B}{n+1})$ and $\Delta^+(\GR{C}{n+1})$.
\end{proof*}%
\begin{ex}
Consider two modifications of $\Delta^+$.

{\sf 1.} Replace $\Delta^+$ with $\widetilde\Delta^+=
\Delta^+\cup\{0\}$, where $\{0\}$ is regarded as the unique
minimal element in this new poset. Hence $\gH(\widetilde\Delta^+)$ 
gains $n$ new edges connecting $\{0\}$ with the simple roots. 
Therefore
\begin{gather*}
   \cku_{\widetilde\Delta^+}(q)=\cku_{\Delta^+}(q)+n(q-1)+1\\
   \ckl_{\widetilde\Delta^+}(q)=\ckl_{\Delta^+}(q)+q^n \ .
\end{gather*}
It follows that 
\[
  \gD_{\widetilde\Delta^+}(q)=\gD_{\Delta^+}(q)-\frac{q^n-nq+n-1}{(q-1)^2}=
-(q^{n-2}+2q^{n-3}+\ldots +(n-2)q) \ .
\]

{\sf 2.} Assume that $n\ge 2$ and
consider $\Delta^+\setminus\Pi$ as subposet of $\Delta^+$.
Then the minimal elements of $\Delta^+\setminus\Pi$ are the roots of height 2.
Here we obtain
$\cku_{\Delta^+\setminus\Pi}(q)=\cku_{\Delta^+}(q)-(n-1)q^2-1$.
But formulae for $\ckl$ depends on the presence of a branching node in 
the Dynkin diagram, i.e., on the presence of a simple root which is covered
by three roots. More precisely,  
\[
\ckl_{\Delta^+\setminus\Pi}(q)=\ckl_{\Delta^+}(q)-
\left\{\begin{array}{rl} 
(n-2)q^2+2q, &  \mbox{if $\Delta^+$ does not have a branching node}; \\
q^3+(n-4)q^2+3q, &  \mbox{if $\Delta^+$ has a branching node\ .}  
\end{array}\right.
\]
Then
\[
  \gD_{\Delta^+\setminus\Pi}(q)=\left\{\begin{array}{rl}
n-2, &  \mbox{if $\Delta^+$ does not have a branching node}; \\
q+n-2, &  \mbox{if $\Delta^+$ has a branching node\ .}
\end{array}\right.
\]
Thus, the deviation polynomial of $\Delta^+, \Delta^+\cup\{0\}$, and
$\Delta^+\setminus\Pi$ always has the nonzero coefficients of the same sign.
\end{ex}
%\vskip-1ex

%%%%%%%%%%%%%%%%   Section 3

\section{Covering polynomials for the poset of {\sf ad}-nilpotent ideals}
\label{tri}

\noindent
Let $\g$ be the simple complex Lie algebra corresponding to $\Delta$ (if $\Delta$ is reduced).
Fix a triangular decomposition $\g=\ut \oplus\te\oplus\ut^-$, where $\te$
is a Cartan subalgebra and %%%$\ut=\underset{\gamma\in\Delta^+}{\bigoplus}\g_\gamma$.
the set of $\te$-roots in $\ut$ is $\Delta^+$.
Then $\be=\te\oplus\ut$ is the fixed Borel subalgebra.

An {\it \adn ideal} of $\be$ is a subspace $\ce\subset \ut$ such that
$[\be,\ce]\subset\ce$.  Then $\ce$ is a sum of certain root spaces in $\ut$,  
$\ce=\bigoplus_{\gamma\in I}\g_\gamma$.  
Here $I$ is necessarily an upper ideal of $\Delta^+$, and this shows that
the poset of {\sf ad}-nilpotent ideals of $\be$ is isomorphic to
the poset of  {upper} ideals of $\Delta^+$.
It will be denoted by $\AD$ or $\AD(\g)$. If $I\in \AD$ is considered
as a  subset of $\Delta^+$, then $\kappa(I)=\#\min(I)$ and $\iota(I)=
\#\max(\Delta^+\setminus I)$. The elements of $\min(I)$ are called {\it generators\/} 
of $I$. 
For $\gamma\in \max(\Delta^+\setminus I)$, 
%%then $I\cup\{\gamma\}$ is again an upper ideal and 
the passage $I\mapsto I\cup\{\gamma\}$ is called an 
{\it extension\/} of $I$. Thus, $\kappa(I)$ (resp. $\iota(I)$) is the number
of generators (resp. extensions) of $I$.

Recall that $\AD\simeq J^*(\Delta^+)$ and $\ck_\AD$ is the covering polynomial of $\AD$.
By Theorem~\ref{ravno1}, $[q^k]\ck_\AD$  equals
the number of $k$-element antichains in $\Delta^+$.
%We write $\ck_\AD$ for this polynomial. 
Here $\deg\ck_\AD=\rk\Delta=n$.
Explicit formulae for $\ck_{\AD}$ for all simple Lie algebras $\g$ can be found e.g. in 
\cite[Section\,6]{duality}.
The  polynomials $\ck_{\AD}$ occur in various contexts, see
\cite{ath1,bessis,duality,charney}.
There is a uniform expression 
(and proof) for the number of all {\sf ad}-nilpotent ideals,
i.e., $\ck_\AD(1)$, see \cite{cp2}. 
Since $\ck_\AD$ is palindromic,
$\ck'_\AD(1)=\frac{n}{2}\ck_\AD(1)$, which yields the expression for the
number of edges in $\gH(\AD)$, see \cite{rodstv}. The coefficients of $\ck_\AD$ are of great interest; for instance, for $\Delta$ of type $\GR{A}{n}$, one obtains the classical 
{\it Narayana numbers}.
But no uniform approach to describing the coefficients of $\ck_\AD$ is known.
For future use, we record the relation between the number of
vertices and edges in $\gH(\AD)$:
\begin{equation}  \label{drob-AD}
\frac{\#\EE(\AD)}{\#\AD}=\frac{\ck'_\AD(1)}{\ck_\AD(1)}=\frac{n}{2}=
\frac{\#(\Delta^+)}{h}\ .
\end{equation}
%%
%%If $\PP=\AD(\g)$, then $\kappa(\ce)$ is the number of generators of $\ce$
%%for any \adn ideal $\ce$. So, in this case, the covering polynomial of 
%%$\AD(\g)$ is the generalised Narayana polynomial $\N_\g$
%%from Section~\ref{pos_all}. These polynomials are written out in 
%%\cite[Section\,6]{duality}.
%%
There is no Lie algebra associated with the root system $\GR{BC}{n}$, but  
one can still consider the poset of upper ideals in 
$\Delta^+(\GR{BC}{n})$, denoted by $\AD(\GR{BC}{n})$.
By Lemma~\ref{poset-BC},  $\#(\AD(\GR{BC}{n}))$
equals to the number of upper ideals 
in $\Delta^+(\GR{B}{n+1})\setminus \Pi$. The latter is known to be equal
to $\genfrac{(}{)}{0pt}{}{2n+1}{n}$ \cite{eric}.
%%%However, the covering polynomial of $\AD(\GR{BC}{n})$ did not appear in print.

\begin{prop}   \label{cover-BC}
The covering polynomial of\/ $\AD(\GR{BC}{n})$ equals
$\sum_{k\ge 0} \genfrac{(}{)}{0pt}{}{n}{k}\genfrac{(}{)}{0pt}{}{n+1}{k}
q^k$.
\end{prop}
\begin{proof} The poset $\Delta^+(\GR{BC}{n})$ is 
isomorphic to the trapezoidal poset $\mathcal T(n,n+1)$ considered by Stembridge
\cite{stembr}. Therefore,  the coefficient of
$q^k$ in the covering polynomial of $\AD(\GR{BC}{n})$ equals the number of
$k$-element antichains in $\mathcal T(n,n+1)$.
By \cite[Theorem~5.4]{stembr}, the latter is the same as the number
of $k$-element antichains in the rectangular poset $\mathcal R(n,n+1)$.
It is easily seen that the number of $k$-element antichains in 
$\mathcal R(n,n+1)$ is $\genfrac{(}{)}{0pt}{}{n}{k}\genfrac{(}{)}{0pt}{}{n+1}{k}$.
\end{proof}

\noindent
Unlike the covering polynomial for the poset of upper ideals in a 
reduced root system, this polynomial is not palindromic.

\begin{cl}    \label{edges:BC}
The number of edges of\/ $\gH(\AD(\GR{BC}{n}))$ is equal to 
$(n+1)\genfrac{(}{)}{0pt}{}{2n}{n+1}=n\genfrac{(}{)}{0pt}{}{2n}{n}$.
\end{cl}
\begin{proof}
%%\noindent
$\displaystyle  \frac{d}{dq}
\bigl(\sum_{k\ge 0} \genfrac{(}{)}{0pt}{}{n}{k}\genfrac{(}{)}{0pt}{}{n+1}{k}q^k
\bigr)
\vert_{q=1}=\sum_{k\ge 0}k\genfrac{(}{)}{0pt}{}{n}{k}\genfrac{(}{)}{0pt}{}{n+1}{k}
=(n+1)\sum_{k\ge 0}\genfrac{(}{)}{0pt}{}{n}{k-1}\genfrac{(}{)}{0pt}{}{n}{k}=\\ \hfil
(n+1)\genfrac{(}{)}{0pt}{}{2n}{n-1}$.
\end{proof}%

\noindent
The poset $\AD(\GR{BC}{n})$ can be regarded as a particular case of the following 
series of examples.
An \adn ideal $\ce$ is said to be {\it strictly positive\/}, if $\ce\subset [\ut,\ut]$. The combinatorial counterpart is that an
upper ideal $I\subset \Delta^+$ is {strictly positive}, 
if $I\cap \Pi=\varnothing$. The corresponding sub-poset of $\AD$ is denoted by
$\AD_0$ or $\AD_0(\g)$.
Clearly, $\AD_0\simeq J^*(\Delta^+\setminus \Pi)$ is a distributive lattice.
%whose poset of meet-irreducible elements is isomorphic to $\Delta^+\setminus \Pi$. 

By Lemma~\ref{poset-BC}, $\AD(\GR{BC}{n})\simeq \AD_0(\GR{B}{n+1}) 
\simeq \AD_0(\GR{C}{n+1})$. This prompts a natural question about ${\AD_0}(\g)$ for
the other root systems (simple Lie algebras). 
A uniform formula for $\#\AD_0(\g)$, i.e., for $\ck_{\AD_0(\g)}(1)$,
is found by Sommers~\cite{eric}. 
In our setting, we are interested in the covering polynomial
$\ck_{\AD_0(\g)}$. % and the number of edges of $\gH(\AD_0(\g))$\,?
The answer for $\AD_0(\GR{B}{n+1})$ and $\AD_0(\GR{C}{n+1})$ 
is given in  
%Lemma~\ref{poset-BC},
Proposition~\ref{cover-BC}. %, and Corollary~\ref{edges:BC}.
The case of $\g=\slno$ is easy, because  $\Delta^+(\GR{A}{n})\setminus\Pi
\simeq \Delta^+(\GR{A}{n-1})$ and
hence $\ck_{\AD_0(\GR{A}{n})}=\ck_{\AD(\GR{A}{n-1})}$. 
The exceptional root systems %(simple Lie algebras) 
can be handled directly. 
For $\GR{D}{n}$, %and $\GR{E}{n}$, 
the answer is not easy to obtain. 

Conjecturally, $\ck_{\AD_0}$ can be expressed via Chapoton's $\eus F$-triangle as follows.
Let $C_{k,l}$ be the set of cones in the cluster complex of $\Delta$ spanned by
$k$ positive roots and $l$ negative simple roots, and $f_{k,l}=\# (C_{k,l})$.
Here $f_{k,l}=0$ if $k+l>n$. Define the  $\eus F$-triangle by its generating function
\[
   \eus F(\Delta)=\eus F(x,y)=\sum_{k,l\ge 0} f_{k,l}x^k y^l .
\]
We refer to  \cite{chap04}  for relevant definitions and other background.
Then we conjecture that
\[  %beq   \label{eq:F1}
    \ck_{\AD_0}(q)= \sum_{k\ge 0} f_{k,0}\,q^k (1-q)^{n-k}=(1-q)^n \eus F(\frac{q}{1-q},0)    
\]
or 
\[ %beq  \label{eq:F2}
    \ck_{\AD_0}(q)= \sum_{k,l\ge 0} f_{k,l}(-1)^l (q-1)^{n-k}=(q-1)^n \eus F(\frac{1}{q-1},-1) .   
\]
This gives the correct formula for $\ck_{\AD_0}(q)$ whenver we can verify it.
(Explicit formulae for the $\eus F$-triangle can be found  in \cite{chap02} or \cite{kratt}.)
We also notice that the coefficients of $\ck_{\AD_0}(q)$ give the "very positive $H$-vector"
in Chapoton's terminology in \cite{chap02}.
%This holds to be true for $\GR{A}{n}$, $\GR{B}{n}$, $\GR{C}{n}$,
%$\GR{F}{4}$, $\GR{G}{2}$. It is also known that the usual ``$h$-vector"
%gives the coefficients of $\ck_{\AD}(q)$.
For $\GR{D}{n}$, this yields  the following conjectural expression:
\[
   [q^k]\ck_{     \AD_0(\GR{D}{n})         }=\bigl(\genfrac{(}{)}{0pt}{}{n-1}{k}\bigr)^2+
   \frac{k-2}{n-1}\genfrac{(}{)}{0pt}{}{n-1}{k} \genfrac{(}{)}{0pt}{}{n-1}{k-1} .
%\frac{1}{n}\genfrac{(}{)}{0pt}{}{n}{k}\bigl( (k+1)\genfrac{(}{)}{0pt}{}{n-2}{k+1}+
%(k+2)\genfrac{(}{)}{0pt}{}{n-2}{k}+ (k-1)\genfrac{(}{)}{0pt}{}{n-2}{k-1} \bigr).
\]
The information for the exceptional Lie algebras is gathered in Table~\ref{table:AD0}.
\begin{table}[htb]  
\begin{tabular}{cl}
$\Delta$ & \qquad\qquad $\ck_{\AD_0}(q)$  \\ \hline
$\GR{E}{6}$ & $1+\phantom{1} 30q+\phantom{1} 135q^2+\phantom{1} 175q^3+
\phantom{15} 70 q^4+\phantom{113} 7 q^5$ 
{\rule{0pt}{2.4ex}} 
\\
$\GR{E}{7}$ & $1+\phantom{1} 56q+\phantom{1} 420q^2+\phantom{1} 952q^3+
\phantom{1} 770 q^4 +\phantom{1} 216 q^5 +\phantom{1} 16 q^6$ 
{\rule{0pt}{2.4ex}} \\
$\GR{E}{8}$ & $1+112q+1323 q^2+4774q^3 +6622 q^4+3696 q^5+770q^6+44q^7$ 
\\
$\GR{F}{4}$ & $1+\phantom{1} 20q+\phantom{51} 35q^2+\phantom{15}10 q^3$ 
\\
$\GR{G}{2}$ & $1+\phantom{11} 4q$ 
\\ 
\hline
\end{tabular}
\vskip.7ex 
\caption{ The covering polynomials for $\AD_0(\g)$, $\g$ being exceptional}
\label{table:AD0}
\end{table}%

\noindent
Using these data and above information for the classical series,
%%of Table~\ref{table:AD0}, 
we obtain:
\begin{equation}  \label{drob-AD0}
\frac{\#\EE(\AD_0)}{\#(\AD_0)}=\frac{\ck'_{\AD_0}(1)}{\ck_{\AD_0}(1)}=
\frac{n}{2}\cdot \frac{h-2}{h-1}=
\frac{\#(\Delta^+\setminus\Pi)}{h-1}\ .
\end{equation}
This equality has a striking similarity with~\eqref{drob-AD}, and it would be interesting to find
a conceptual explanation for it.  More precisely, one can suggest a general pattern behind
\eqref{drob-AD} and \eqref{drob-AD0}.
Let $\LL$ be a graded poset such that the maximal length of chains in $\LL$ is $r-1$
(i.e., such a chain contains $r$ elements). These two equalities are manifestations of the following phenomenon:

{\it For some ``good'' graded posets $\LL$, one has} \quad
$\displaystyle       
    \frac{\# \EE(J^*(\LL))}{\# J^*(\LL)}=\frac{\# \LL}{r+1}$.

\noindent Besides $\Delta^+$ and $\Delta^+\setminus \Pi$,  the weight posets of some
(finite-dimensional) representations of $\g$ considered in \cite{mmj}
also have this property.

%%%%%%%%%%%%%%%%   Section 4
\section{Covering polynomials for the poset of abelian ideals}
\label{genera}

\noindent
There is an interesting subposet of $\AD$, where the two covering 
polynomials are different.
An \adn ideal $\ce\subset \be$ is {\it abelian\/} if $[\ce,\ce]=0$. The combinatorial counterpart is
that an upper ideal $I\subset\Delta^+$ is said to be {\it abelian}, if 
$\gamma'+\gamma''\not\in \Delta^+$ for each pair $\gamma',\gamma''\in I$.
Let $\Ab=\Ab(\g)$ be the subposet of 
$\AD$ consisting of all abelian ideals.
Clearly $\Ab$ is a graded meet-semilattice. It follows that $\Ab$
is a (distributive) lattice if and only if there is a unique maximal abelian ideal,
which happens only for $\GR{C}{n}$ and $\GR{G}{2}$.
In all other cases the upper and lower covering polynomials are different.

Our next results rely on
the relationship, due to D.\,Peterson, between the abelian ideals and the so-called
{\it minuscule elements\/} of the affine Weyl group of $\g$. Recall the
necessary setup.
\\[.6ex]
We have the real vector space $V=\oplus_{i=1}^n{\mathbb R}\ap_i$ and 
$(\ ,\ )$ a $W$-invariant inner product on $V$. 

$Q=\oplus _{i=1}^n {\mathbb Z}\ap_i \subset V$ is the root lattice; 

$Q^+=\{\sum_{i=1}^n m_i\ap_i \mid m_i=0,1,2,\dots \}$ is the monoid generated by the positive
roots.

\noindent
As usual, $\mu^\vee=2\mu/(\mu,\mu)$ is the coroot
for $\mu\in \Delta$ and $Q^\vee=\oplus _{i=1}^n {\mathbb Z}\ap_i^\vee$  
is the coroot lattice. % in $V$.
\\[.6ex]
Letting $\widehat V=V\oplus {\mathbb R}\delta\oplus {\mathbb R}\lb$, we extend
the inner product $(\ ,\ )$ on $\widehat V$ so that $(\delta,V)=(\lb,V)=
(\delta,\delta)=
(\lb,\lb)=0$ and $(\delta,\lb)=1$.
Then 
\begin{itemize}
\item[] \ 
$\widehat\Delta=\{\Delta+k\delta \mid k\in {\mathbb Z}\}$ is the set of affine
(real) roots; 
\item[] \ $\HD^+= \Delta^+ \cup \{ \Delta +k\delta \mid k\ge 1\}$ is
the set of positive affine roots; 
\item[] \ $\widehat \Pi=\Pi\cup\{\ap_0\}$ is the corresponding set
of affine simple roots. 
\end{itemize}
Here  $\ap_0=\delta-\theta$.
For $\ap_i$ ($0\le i\le n$), let $s_i$ denote the corresponding 
reflection in $GL(\HV)$.
That is, $s_i(x)=x- (x,\ap_i)\ap_i^\vee$ for any $x\in \HV$.
The affine Weyl group, $\HW$, is the subgroup of $GL(\HV)$
generated by the reflections $s_i$, $i=0,1,\ldots,n$.
If the index of $\ap\in\widehat\Pi$ is not specified, then we merely write
$s_\ap$. 
The inner product $(\ ,\ )$ on $\widehat V$ is
$\widehat W$-invariant. The notation $\beta>0$ (resp. $\beta <0$)
is a shorthand for $\beta\in\HD^+$ (resp. $\beta\in -\HD^+$).
The length function on $\widehat W$ with respect
to  $s_0,s_1,\dots,s_p$ is denoted by $\ell$.
For $w\in\HW$, we set $\eus N(w)=\{\nu\in\HD^+\mid w(\nu)<0\}$.
Then $\# \eus N(w)=\ell(w)$.

\begin{df}[D.\,Peterson]
An element $w\in \HW$ is said to be {\it minuscule\/}, if 
$\eus N(w)=\{\delta-\gamma\mid \gamma\in I\}$ for some subset $I\subset \Delta$.
\end{df}%
Then one can easily show that $I\subset \Delta^+$, $I$ is an abelian ideal, and
this correspondence yields a bijection between the minuscule elements of
$\HW$ and the abelian ideals. 
Furthermore, if $w$ is minuscule and $w^{-1}(\ap)=-\mu+k\delta$
($\ap\in\HP$, $\mu\in \Delta$), then $k\ge -1$.
(More generally, this holds for elements of $\HW$ corresponding to arbitrary
\adn ideals, see \cite{cp1}).
If $w$ is minuscule, then $I_w$ denotes the corresponding abelian ideal.
Conversely, given $I\in\Ab$, then $w_I$ stands for the corresponding minuscule
element. If $I\in\Ab$ and $\gamma\in \max(\Delta^+\setminus I)$, then the ideal
$I'=I\cup\{\gamma\}$ is not necessarily abelian. In this section, we are interested
only in abelian extensions, i.e., those with abelian $I'$.

\begin{lm}  \label{theta2}
Suppose $w\in\HW$ is minuscule and $w^{-1}(\ap)=-\mu+2\delta$, where 
$\ap\in\HP$ and $\mu\in\Delta$. Then $\mu=\theta$.
\end{lm}
\begin{proof}
We have $\ap=w(2\delta-\mu)=w(2\delta-\theta)+w(\theta-\mu)$. Here
$\theta-\mu\in Q^+$. Hence both $2\delta-\theta$ and $\theta-\mu$ do not
belong to $\eus N(w)$. Therefore if $\theta\ne\mu$, then one obtains a contradiction
with the fact that $\ap$ is simple.
\end{proof}

\noindent
For a minuscule $w$, consider the vector $\mathbf k=\mathbf k_w=
(k_0,k_1,\ldots,k_n)$, 
where $k_i$ is defined by the equality $w^{-1}(\ap_i)=-\mu_i+k_i\delta$
($\mu_i\in\Delta$). Recall that $k_i\ge -1$ for all $i$. 

\begin{prop}    \label{ogr}   %\leavevmode\par
We have
\begin{itemize}
\item[\sf (i)] \ $k_i\le 2$ for all $i$;
\item[\sf (ii)] \ there is at most one index $i$ such that $k_i=2$.
The corresponding simple root $\ap_i$ is necessarily long.
\item[\sf (iii)] \ $k_0\le 1$; that is, $k_0\ne 2$.
\end{itemize}
\end{prop}
\begin{proof}
(i) \ If $w^{-1}(\ap_i)=-\mu_i+k_i\delta$ and $k_i\ge 3$, then
$w(2\delta-\mu_i)=-(k_i-2)\delta+\ap_i <0$. Hence $w$ is not minuscule.

(ii) \ If $w^{-1}(\ap_i)=-\mu_i+2\delta$, then $\mu_i=\theta$ by Lemma~\ref{theta2}.
Hence such $i$ is unique and $\|\ap_i\|=\|\theta\|$, i.e., $\ap_i$ is long.

(iii) \ Suppose $w^{-1}(\ap_0)=-\mu_0+2\delta$. Then $\mu_0=\theta$
and $w(2\delta-\theta)=\delta-\theta$. However, it is shown in
\cite[Prop.\,2.5]{imrn} that $w(2\delta-\theta)\in\Delta^+$ for any non-trivial
minuscule element $w$.
\end{proof}

\noindent
We shall say that $\mathbf k_w$ is the {\it shift vector\/} of $w$
or $I_w$. If $w=w_I$ ($I\in\Ab$), then
we also write $\mathbf k_I$ for this vector. 

\begin{thm}  \label{ext}
Let $(k_0,k_1,\ldots,k_n)$ be the shift vector of $I\in\Ab$. Then 
$\kappa(I)=\#\{i \mid k_i=-1\}$ and
$\iota(I)=\#\{i \mid k_i=1\}$.
\end{thm}
\begin{proof}
1. It is shown in \cite[Theorem\,2.2]{imrn} that $\gamma\in I$ is a generator
if and only if $w_I(\delta-\gamma)=-\ap_i\in -\HP$. That is, $k_i=-1$ for the
corresponding coordinate $i$.

2. Suppose that  $k_i=1$, i.e., $w_I^{-1}(\ap_i)=-\mu_i+\delta$. 
Then  $w_I(-\mu_i)=\ap_i-\delta <0$. As $w_I$ is minuscule, $\mu_i$ must be positive.
%then one easily sees that $\mu_i\in \Delta^+$. 
Next, $\eus N(s_iw_I)=\eus N(w_I)\cup\{ w_I^{-1}(\ap_i)\}$. Hence
 $s_iw_I$ is again minuscule, and the corresponding abelian
ideal is $\tilde I=I\cup\{\mu_i\}$. Conversely, if $I\to I\cup \{\gamma\}$ is an
abelian extension of $I$, then $w_{\tilde I}=s_iw_I$ for 
some $i\in \{0,1,\ldots,n\}$ and 
$w_{\tilde I}(\delta-\gamma)=-\ap_i$ \cite[Theorem\,2.4]{imrn}. Then 
$w_I(\delta-\gamma)=\ap_i$, i.e., $k_i=1$.
\end{proof}

\noindent
Since the minuscule elements ($\sim$ abelian ideals) can be constructed 
recursively, we obtain, as a consequence of this theorem,
a method to compute recursively  the shift vector.
One starts with the minuscule element $1\in\HW$ (or the empty ideal).
The corresponding shift vector is $\mathbf k_{1} =(1,0,\ldots,0)$. The inductive step is to replace $w=w_I$ with $s_iw$ for some $i$.
However one has to be careful while choosing  $s_i$, otherwise $s_iw$ may fail to
be minuscule.

\begin{prop}   \label{ind_proc}
Suppose that $w\in\HW$ is minuscule. Then

\textbullet \quad  $s_jw$ is again minuscule if and only if\/
$(\mathbf k_w)_j=1$. 

\textbullet \quad If\/ $(\mathbf k_w)_j=1$, then 
$\mathbf k_{s_jw}=\mathbf k_w-
(\text{the $j$-th column of the extended Cartan matrix of $\Delta$} )^t$.
\end{prop}
\begin{proof}
The first claim is essentially proved in the second part of the above theorem.
The second claim follows from the assumption $(\mathbf k_w)_j=1$ and 
the equalities:

\centerline{
$(s_jw)^{-1}(\ap_i)=w^{-1}(\ap_i)-(\ap_i,\ap_j^\vee)w^{-1}(\ap_j); \quad i=0,1,\ldots,n$.
}

\noindent
Recall that the extended Cartan matrix is the $(n{+}1)\times (n{+}1)$ matrix
with  entries $(\ap_i,\ap_j^\vee)$, $0\le i,j\le n$.
\end{proof}

\noindent
{\bf Remark.} Write $\theta=\sum_{i =1}^n c_i\ap_i$ and set  $c_0=1$.
Then $\sum_{i=0}^n c_i\ap_i=\delta$. Since $\delta$ is $\HW$-invariant, 
the definition of $k_i$'s implies that $\sum_{i=0}^n c_ik_i=1$ for any
shift vector. Hence $\mathbf k_I$ is fully determined by $k_1,\ldots,k_n$.
Let $z=z_I\in V$ be the unique point such that $(z_I,\ap_i)=k_i$,
$i=1,\ldots,n$. Then  $z\in Q^\vee$. (Again, this is true in the context 
of arbitrary \adn ideals, see \cite{cp2}). Note that $k_0=1-(z,\theta)$.
The constraints of Proposition~\ref{ogr} show that 
$-1\le (z,\ap_i)\le 2$ for $i=1,\ldots,n$ and $0\le (z,\theta)\le 2$.
Furthermore, a stronger result is valid.
It was shown by Kostant \cite{ko1} that the mapping
$I\in\Ab \mapsto z_I\in V$ sets up a bijection between the abelian ideals
and the points $z\in Q^\vee$ such that $-1\le (z,\gamma)\le 2$
for all $\gamma\in \Delta^+$.

Let $\Delta^+_l$ denote the set of long positive roots
and $\Pi_l:=\Delta^+_l\cap \Pi$.
In \cite{imrn}, we constructed a disjoint partition 
of $\Abo:=\Ab\setminus \{\varnothing\}$ parametrised by $\Delta^+_l$.
%%%the long positive roots.
In other words, there is a natural surjective mapping
$\tau:\Abo\to \Delta^+_l$. Given $I\in\Abo$ and the corresponding
minuscule element $w\in \HW$, we set $\tau(I)=w(2\delta-\theta)$.
By \cite[Prop.\,2.5]{imrn}, it is an element of $\Delta^+_l$.
Then $\Ab_\mu=\tau^{-1}(\mu)$.

\begin{rema} Using the above definition of the shift vector of an abelian ideal
and Lemma~\ref{theta2}, one observes that $(\mathbf k_I)_i=2$ if and only if 
$\tau(I)=\ap_i$, i.e., $I\in\Ab_{\ap_i}$.
\end{rema}

One of the main results of  \cite{imrn} is that each $\Ab_\mu$ has a unique maximal
and unique minimal element, and that the maximal elements of $\Ab$ are exactly
the maximal elements of $\Ab_\ap$, $\ap\in\Pi_l$. 
Let $I\in\Abo$ and $\tau(I)=\mu$. By \cite[Prop.\,3.2]{imrn}, if $I\to I'$ is an
extension, then $\tau(I')\curle \tau(I)$. Hence $I$ has an extension outside 
$\Ab_\mu$ only if $\mu\not\in\Pi$.
Now we make that analysis more precise by showing that the number of possible 
abelian extensions of $I\to I'$ such that $I'\not\in\Ab_\mu$  
depends only on $\mu$ and not on $I$.

\begin{thm}  \label{main_lower}
For any $\mu\in\Delta^+_l\setminus \Pi$ and any $I\in\Ab_\mu$, the
%%there is one and the same
number of abelian extensions $I\to I'$ such that $I'\not\in\Ab_\mu$
equals the number of $\ap\in\Pi$ such that $(\mu,\ap)>0$.
%%the number of roots of the form $\gamma=\mu-\ap$ with $\ap\in\Pi_l$.
In particular, this number does not depend on $I$, and if $\Delta$ is simply laced,
then it is equal to $\kappa(\mu)$.
\end{thm}
\begin{proof*}
1. Suppose $(\ap,\mu)>0$ and $\mu':=s_\ap(\mu)=\mu-n_\ap\ap$. Here $n_\ap\in\mathbb N$
and $n_\ap=1$
if and only if $\ap$ is long. Anyway, $\mu'':=\mu-\ap$ is again a root,
and we will work with the sum $\mu=\ap+\mu''$.
Then $w_I^{-1}(\mu''+\ap)=2\delta-\theta$ and,
as shown in the proof of Theorem~2.6
in \cite{imrn}, one has $w_I^{-1}(\mu'')=\delta-\gamma''$ and 
$w_I^{-1}(\ap)=\delta-\gamma'$ for some $\gamma',\gamma''$ such that
$\gamma'+\gamma''=\theta$. This shows that $w=s_\ap w_I$ is minuscule,
$I_w=I\cup\{\gamma'\}$, and  $\tau(I_w)=s_\ap w_I(2\delta-\theta)=s_\ap(\mu)=\mu'$.

2. Conversely, suppose $I\in\Ab_\mu$ and
$I\to I'$ is an extension. Then $I'=I\cup\{\gamma\}$ for
some $\gamma\in\Delta^+$ and $w_{I'}=s_\ap w_I$, where $\ap\in\Pi$ is determined by
the equality $w_I(\delta-\gamma)=\ap$. The condition $\tau(I')\ne \mu$ 
means $s_\ap(\mu)\ne \mu$, i.e., $(\ap,\mu)\ne 0$. To compute the sign, we notice
that $(\ap,\mu)=(w_I^{-1}(\ap), w_I^{-1}(\mu))=(\delta-\gamma, 2\delta-\theta)=
(\gamma, \theta)$, which cannot be negative.
%% and $I'\in \Ab_{\mu'}$ with $\mu\ne\mu'$,
%%then $w_{I'}=s_\ap w_I$ for some $\ap\in\Pi$ and $\tau(I')=s_\ap(\mu)$.
\end{proof*}%

\noindent
\begin{cl}   \label{ext-1}
The ideals having a unique abelian extension are the following:

(a) \ $\varnothing$;

(b)  The maximal elements of posets $\Ab_\mu$, where $\mu\not\in\Pi$ 
and the inequality $(\ap,\mu)>0$ holds for a unique $\ap\in\Pi$;

(c)  The ideals having a unique abelian extension inside $\Ab_\ap$,
$\ap\in\Pi_l$.
\end{cl}%
It was noticed in \cite{imrn} that each $\Ab_\mu$ is a {\it minuscule poset},
i.e., there is a simple Lie algebra $\el$ and a parabolic subalgebra
$\p\subset\el$ with abelian nilpotent radical $\p^{nil}$ such that
$\Ab_\mu$ is isomorphic to the poset of abelian $\p$-ideals in $\p^{nil}$.
The construction of $\el$ as a subalgebra of $\g$ is given in \cite[Section\,5]{imrn}.
Since the structure of the minuscule posets is well known, Theorem~\ref{main_lower}
provides an effective tool for computing the lower covering polynomial in the small rank
cases, e.g. for the exceptional Lie algebras.

\noindent
We say that a root $\gamma\in\Delta^+$ is {\it commutative}, if the upper
ideal generated by $\gamma$ is abelian. Clearly, the set of commutative roots
forms an upper ideal. A uniform description of this ideal (and its cardinality)
is given in \cite[Theorem\,4.4]{rodstv}. 
In particular, if the Dynkin diagram has no branching node
then the number of commutative roots is $n(n+1)/2$.
A trivial but useful remark is that any abelian ideal consists of commutative roots.
\\[.6ex]
In the following assertion, we gather some information that is helpful in
practical computations of the covering polynomials.

\begin{prop}   \label{practical}
\leavevmode\par
\begin{itemize}
\item[\sf (i)] \  $[q^0]\cku_\Ab=1$, $[q^0]\ckl_\Ab=\#(\Pi_l)$;
\item[\sf (ii)] \ $[q]\cku_\Ab=\text{the number of commutative roots}$; 
\item[\sf (iii)] \ $\deg \cku_\Ab \le \deg\ckl_\Ab$. If these degrees are equal, 
           say to $m$, then $[q^m]\cku_\Ab \le [q^m]\ckl_\Ab$;
\item[\sf (iv)] \  $\cku_\Ab(1)=\ckl_\Ab(1)=2^n$;
\item[\sf (v)] \  $\cku'_\Ab(1)=\ckl'_\Ab(1)=(n+1)2^{n-2}$.
\end{itemize}
\end{prop}
\begin{proof}
(i) \ These are the numbers of minimal and maximal elements of $\Ab$.
\\
(ii) \ Obvious.
\\
(iii) \ Let $I$ be an abelian ideal with $m$ generators, say $\{\gamma_1,\ldots,
\gamma_m\}$. Then $I\setminus \{\gamma_1,\ldots,\gamma_m\}$ has at least $m$ 
extensions. Then take $m=\deg\cku_\Ab$.
\\
(iv) \ This is Peterson's theorem on $\#(\Ab)$, see e.g. \cite[Theorem\,2.9]{cp1};
\\
(v) \ This is the number of edges of $\gH(\Ab)$, which is computed in 
\cite[Theorem\,4.1]{rodstv}.
\end{proof}%

\noindent
Let $m$ be the maximal size of an antichain of commutative roots.
Then $\deg\ckl_\Ab\le m$. But this bound may not be sharp, since the ideal
having a prescribed antichain as the set of either generators or minimal elements 
in the complement can be non-abelian. However, in case $\GR{E}{7}$ and
$\GR{E}{8}$ this bound does give the exact value for both degrees, and we obtain
$\deg\ck_{\Ab(\GR{E}{7})}=\deg\ck_{\Ab(\GR{E}{8})}=4$, where $\ck$ is either
$\cku$ or $\ckl$.

Also, it follows from a result of Sommers \cite[Theorem\,6.4]{eric} that 
the number of generators of an abelian ideal
is at most the maximal number of pairwise orthogonal roots in $\Pi$.
This provides an upper bound for $\deg\cku_\Ab$. 
Since it is easy to
find an abelian ideal with such a number of generators, one actually obtains
the exact value of the degree. For instance, $\deg\cku_{\Ab(\GR{D}{n})}=
\left[\frac{n}{2}\right]+1$ and $\deg\cku_{\Ab(\GR{E}{6})}=3$.
If $\deg\ckl_\Ab \le 3$, then both covering polynomials can be computed using
Proposition~\ref{practical}.
This applies to $\GR{F}{4}$ and
$\GR{E}{6}$. For $\GR{E}{7}$ and $\GR{E}{8}$, it suffices to determine one more
value (or coefficient) of $\cku_\Ab$ and $\ckl_\Ab$.

\begin{ex}
$\g=\GR{E}{8}$. 

(1) \ We use Corollary~\ref{ext-1} to compute the coefficient 
$[q]\ckl_{\Ab(\GR{E}{8})}$.
Here the number of roots $\mu$ such that $\kappa(\mu)=1$ is 21
(see Section~\ref{dva}).  The posets $\Ab_{\ap_i}$, $\ap_i\in\Pi$, have the
cardinalities $1,\,2,\,3,\,4,\,5,\,6,\,8,\,6$. Furthermore, for $\GR{E}{8}$, each poset
$\Ab_\mu$ is totally ordered. Hence the contribution from part (c) of the
corollary is $0+1+2+3+4+5+7+5=27$. Thus, $[q]\ckl_{\Ab(\GR{E}{8})}
=1+21+27=49$.

(2) \ Since each $\Ab_\mu$ is a chain, any ideal $I$ has at most one extension
inside its own poset $\Ab_\mu$. Because $\kappa(\mu)\le 3$ for all $\mu\in\Delta^+$,
Theorem~\ref{main_lower} shows that 
\[
  [q^4]\ckl_{\Ab(\GR{E}{8})}=\sum_{\kappa(\mu)=3} \#(\Ab_\mu)-1 \ .
\]
We have $\#\{\mu\in\Delta^+\mid \kappa(\mu)=3\}=[q^3]\cku_{\Delta^+(\GR{E}{8})}=21$.
Of these 21 roots, there are

11 roots with $\#(\Ab_\mu)=1$ \ (these are exactly the roots $\mu$ with 
$(\theta,\mu)\ne 0$, see \cite[5.1]{imrn};

5 roots with $\#(\Ab_\mu)=2$;

3 roots with $\#(\Ab_\mu)=3$;

2 roots with $\#(\Ab_\mu)=4$.
\\[.6ex]
Hence $[q^4]\ckl_{\Ab(\GR{E}{8})}=17$.

(3) Then using Proposition~\ref{practical}(iv),(v), we compute
$\ckl_{\Ab(\GR{E}{8})}(q)=8+49q+87q^2+95q^3+17q^4$.
\end{ex}
%\vskip-1.5ex

%%%%%%%%%%%%%%%%   Section 5

\section{Computing the covering polynomials for $\Ab(\g)$, $\g$ being classical}
\label{classical}

\noindent
In this section, we prove four theorems for four classical series of simple Lie algebras
(root systems). Our proofs are based
on an explicit presentation of the upper ideal of commutative roots and 
understanding which ideals inside it are really abelian. To this end, one has to know 
the generators of all maximal abelian ideals, 
which are determined in \cite{R} (see also \cite[Table\,1]{pr}).
				   
\begin{thm}  \label{thm:covA}
If\/ $\g=\slno$, i.e., $\Delta$ is of type $\GR{A}{n}$, then
\begin{itemize}
\item[\sf (i)] \ $\cku_\Ab(q)=\sum_{k\ge 0}\genfrac{(}{)}{0pt}{}{n+1}{2k}q^k$;
\item[\sf (ii)] \ $\ckl_\Ab(q)=\sum_{k\ge 0}
\bigl(\genfrac{(}{)}{0pt}{}{n}{2k+1}+
\genfrac{(}{)}{0pt}{}{n}{2k-2}\bigr)q^k$;
\item[\sf (iii)] \ $\gD_\Ab(q)=-\sum_{k\ge 0}\genfrac{(}{)}{0pt}{}{n-1}{2k+1}q^k$.
\end{itemize}
\end{thm}%%\begin{proof*}
{\sl Proof.}\quad 
The formula for $\cku$ (but not for $\ckl$!)
is implicit in \cite[Section\,3]{pr}. We recall the necessary
setup and then deduce the expressions (i) and (ii). Then part (iii)
is obtained via formal manipulations.

For $\g=\slno$, each positive root is commutative.
The root $\esi_i-\esi_j\in \Delta^+(\slno)$ is identified with the pair
$(i,j)$.
Suppose $\ah\in\Ab(\slno)$ and $\kappa(\ah)=k$.
That is, $\ah$ has $k$ generators (minimal roots).
If $\min(\ah) 
%%These minimal roots are represented with the collection of numbers
=\{(a_1,b_1),\ldots,(a_k,b_k)\}$, where $a_1< a_2<\ldots < a_k$,
then we actually have
$1\le a_1< a_2<\ldots < a_k< b_1<\ldots <b_k\le n+1$.
Thus, any $2k$-element subset of $[1,n+1]$ gives rise to an abelian ideal 
with $\kappa(\ah)=k$ and vice versa. This yields (i).

The edges of $\gH(\Ab)$ originating in $\ah$ bijectively
correspond to the maximal roots $\gamma$ in $\Delta^+\setminus\ah$ such that
$\{\gamma\}\cup\ah$ is again an abelian ideal.
The set of such maximal roots always contains
$\{(a_1+1,b_2-1),\ldots,(a_{k-1}+1,b_k-1)\}$;
furthermore, if $a_k+1<b_1$, then two more roots are admissible:
$(1,b_1-1), (a_k+1,n+1)$.
From this we deduce that $\iota(\ah)=k$ if and only if
one of the following two conditions hold:

$(\Diamond_1)$ \ $\#\min(\ah)=k+1$ and $a_{k+1}+1=b_1$. 

$(\Diamond_2)$ \ $\#\min(\ah)=k-1$ and $a_{k-1}+1<b_1$.
\\[.7ex]
In case $(\Diamond_1)$ the ideal is determined by a sequence of 
$2k+1$ integers

$1\le a_1<\ldots < a_k< a_{k+1}=b_1-1< b_2-1<
\ldots <b_{k+1}-1\le n$. 
\\[.6ex]
Hence there are $\genfrac{(}{)}{0pt}{}{n}{2k+1}$
such possibilities.
\\[.6ex]
In case $(\Diamond_2)$ the ideal is determined by a sequence of 
$2k-2$ integers

$1\le a_1<\ldots < a_{k-1}<b_1-1< b_2-1<
\ldots <b_{k-1}-1\le n$. 
\\[.6ex]
Hence there are $\genfrac{(}{)}{0pt}{}{n}{2k-2}$
such possibilities. This proves (ii).

(iii) \ It follows from parts (i) and (ii) that
\begin{multline*}
\cku_\Ab(q)-\ckl_\Ab(q)= \\  \textstyle
-\sum_{k\ge 0}
\bigl(\genfrac{(}{)}{0pt}{}{n}{2k+1}-\genfrac{(}{)}{0pt}{}{n+1}{2k}
+\genfrac{(}{)}{0pt}{}{n}{2k-2} \bigr)q^k= 
-\sum_{k\ge 0}
\bigl(\genfrac{(}{)}{0pt}{}{n-1}{2k+1}-2\genfrac{(}{)}{0pt}{}{n-1}{2k-1}
+\genfrac{(}{)}{0pt}{}{n-1}{2k-3}\bigr)q^k= \\  \textstyle
-(q-1)^2 \sum_{k\ge 0}\genfrac{(}{)}{0pt}{}{n-1}{2k+1}q^k.\quad \square
\end{multline*}
%\end{proof*}%
% 
\vskip-1ex
\begin{thm}  \label{thm:covC}
If\/ $\g=\spn$, i.e., $\Delta$ is of type $\GR{C}{n}$, then
$\cku_\Ab(q)=\ckl_\Ab(q)=\sum_{k\ge 0}\genfrac{(}{)}{0pt}{}{n+1}{2k}q^k$.
\end{thm}\begin{proof}
The commutative roots are $\{\esi_i+\esi_j\mid 1\le i\le j\le n\}$.
Since these roots form the unique maximal abelian ideal,
the poset $\Ab(\spn)$ is a distributive lattice and
%whose subset of meet-irreducibles is isomorphic to the set of commutative roots. Hence 
$\cku_\Ab(q)=\ckl_\Ab(q)$.
The explicit form of this polynomials stems from the description given in
\cite[Section\,3]{pr}: the ideals with $k$ generators (or with $k$ extensions)
are in a bijection with the the sequences 
$1\le a_1< a_2<\ldots < a_k\le b_1<\ldots <b_k\le n$. Hence there are
$\genfrac{(}{)}{0pt}{}{n+1}{2k}$ possibilities for them.
\end{proof}

\noindent
The poset of commutative roots for $\spn$ is represented by the triangular Ferrers
diagram with row lengths $(n,n-1,\ldots,1)$ (=\,triangle `of size $n$').
In the following two theorems, such a triangle occurs as a subposet of the poset of 
commutative roots
for $\sono$ and $\sone$, which allows us to exploit the
formula of Theorem~\ref{thm:covC}.

\begin{thm}  \label{thm:covB}
If\/ $\g=\sono$, i.e., $\Delta$ is of type $\GR{B}{n}$, then
\begin{itemize}
\item[\sf (i)] \ $\cku_\Ab(q)=\sum_{k\ge 0}\genfrac{(}{)}{0pt}{}{n+1}{2k}q^k$;
\item[\sf (ii)] \ $\ckl_\Ab(q)=\sum_{k\ge 0}
\bigl(\genfrac{(}{)}{0pt}{}{n-1}{2k+1}{+}
\genfrac{(}{)}{0pt}{}{n}{2k-1}{+}\genfrac{(}{)}{0pt}{}{n-1}{2k-2}\bigr)q^k$;
\item[\sf (iii)] \ $\gD_\Ab(q)=-\sum_{k\ge 0}\genfrac{(}{)}{0pt}{}{n-2}{2k+1}q^k$.
\end{itemize}
\end{thm}
\begin{proof}
Here the commutative roots are 
$\{\esi_i+\esi_j\mid 1\le i<j \le n\}\cup
\{\esi_1-\esi_i\mid 2\le i\le n\}\cup\{\esi_1\}$.\\
Graphically, this set is represented by a skew Ferrers diagram
with row lengths $(2n-1,n-2,n-3,\ldots,1)$. See the sample figure for 
%$\g=\mathfrak{so}_{13}$ 
$\Delta$ of type $\GR{B}{6}$,
where $\esi_{ij}=\esi_i-\esi_j$ and $\ov{\esi}_{ij}=\esi_i+\esi_j$.
%the leftmost (resp. rightmost) box is $\ap_1=\esi_1-\esi_2$
%(resp. $\theta=\esi_1+\esi_2$). The central box in the upper row is $\esi_1$,
%and lowest box is $\esi_{n-1}+\esi_n$.

\begin{figure}[htb]
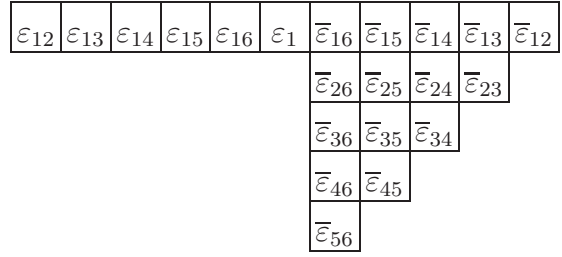

\begin{center}
\setlength{\unitlength}{0.02in}
\Yboxdim3.4ex  
{\young(\zab\zac\zad\zae\zaf\za\zafn\zaen\zadn\zacn\zabn  ,::::::\zbfn\zben\zbdn\zbcn ,::::::\zcfn\zcen\zcdn ,::::::\zdfn\zden ,::::::\zefn )}
\caption{The (po)set of commutative roots 
for $\mathfrak{so}_{13}$} \label{ab:sono}
\end{center}
\end{figure}

\noindent
In such a diagram, the maximal element ($\theta=\ov{\esi}_{12}$) appears in the 
northeast corner and the smaller elements appear to the south and west. The edges 
correspond to the pairs of 
boxes having a common side. The direction of arrows is either `$\to$' or `$\uparrow$'.
This Ferrers diagram consists of the tail of length $n$ and the triangle
`of size $n-1$'. The triangle itself represents an abelian ideal, and 
the structure of the set of ideals sitting inside this triangle
is the same as for $\mathfrak{sp}_{2n-2}$.
\\[.7ex]
(i) \ Let us compute the number of abelian ideals $\ah$ with $k$ generators.

$\bullet$ \ By Theorem~\ref{thm:covC}, 
the number of such ideals inside the triangle is equal to
$\genfrac{(}{)}{0pt}{}{n}{2k}$.

$\bullet$ \ Suppose that $\ah$ has the tail of length $m$, $m\ge 1$, i.e.,
$\ah$ has the generator $\esi_1-\esi_{n+2-m}$ in the upper row.
Then the rest of this row (to the right) is also in the ideal, and
the ideal is determined by its part lying in the triangle of size
$n-2$, in the rows from 2 to $n-1$. The condition of being abelian means that
$\ah$ cannot have elements from $m-1$ leftmost columns of the triangle.
(Formally: if $\esi_1-\esi_{n+2-m}\in \ah$, then for all other roots $\esi_i+\esi_j
\in\ah$, $2\le i\le j$, we must have $j\le n+1-m$.) 
Hence as a degree of freedom for further constructing $\ah$
we have a triangle of size $n-1-m$, where an ideal with $k-1$ generators
has to be chosen. The symplectic case shows that the number of such possibilities
equals $\genfrac{(}{)}{0pt}{}{n-m}{2k-2}$.

Thus, the total number of abelian ideals with $k$ generators equals
\[
    \genfrac{(}{)}{0pt}{}{n}{2k}+\sum_{m\ge 1} \genfrac{(}{)}{0pt}{}{n-m}{2k-2}=
\genfrac{(}{)}{0pt}{}{n}{2k}+\genfrac{(}{)}{0pt}{}{n}{2k-1}=
\genfrac{(}{)}{0pt}{}{n+1}{2k}\ .
\]
(ii) \ Let us compute the number of all abelian ideals $\ah$ with $k$ extensions.
Here the argument is similar in the spirit, but more tedious. 
%%We only sketch the main steps.

$\bullet$ \ Suppose $\ah$ lies in the triangle. The difficulty here is that
$\ah$ may have an extension that does not fit in the triangle
(namely, if the length of the first row equals $n-1$). Therefore the symplectic formula
does not immediately apply.

Let $p\le n-1$ be the length of the first row of $\ah$. Then $\ah$ certainly have the 
extension in the first row. The rest of $\ah$ (in the second row and below) 
sits in the triangle of size $p-1$ and must have $k-1$ extensions.
By the symplectic formula, the number of possibilities here is 
$\genfrac{(}{)}{0pt}{}{p}{2k-2}$. Hence, the total number of possibilities for 
the ideals inside the triangle equals
$\sum_{p\le n-1}\genfrac{(}{)}{0pt}{}{p}{2k-2}=\genfrac{(}{)}{0pt}{}{n}{2k-1}$.
%%The number of such abelian ideals inside the triangle
%%equals $\genfrac{(}{)}{0pt}{}{n}{2k-1}$.

$\bullet$ \ Suppose an abelian ideal $\ah$ has the tail of length $m$, $m\ge 1$.
Let $s$ be the length of the second row of $\ah$.
Then, as explained in the proof of part (i),  $s\le n-m-1$.
Here one has to distinguish two cases. 

(1) \ If $s=n-m-1$, then $\ah$ has no extensions 
in the first two rows.
Hence all $k$ extensions must occur in row  no.\,3 and below. This part
of $\ah$ sits in the triangle of size $n-m-2$. Therefore
one has $\genfrac{(}{)}{0pt}{}{n-m-1}{2k}$
possibilities for constructing an ideal.
    %%$\genfrac{(}{)}{0pt}{}{n-1}{2k+1}$

(2) \ If $s<n-m-1$, then $\ah$ has extensions in both 
the first and second row. Hence the lower part of $\ah$, in row no.\,3 and below,
must have $k-2$ extensions. Since this lower part 
sits inside the triangle of size $s-1$, one has $\genfrac{(}{)}{0pt}{}{s}{2k-4}$
possibilities. Altogether, we obtain 
$\sum_{s\le n-m-2}\genfrac{(}{)}{0pt}{}{s}{2k-4}=\genfrac{(}{)}{0pt}{}{n-m-1}{2k-3}$
variants.
\\[.8ex]
Thus, the total number of abelian ideals with $k$ extensions equals
\[
    \genfrac{(}{)}{0pt}{}{n}{2k-1}+
\sum_{m\ge 1} \genfrac{(}{)}{0pt}{}{n-m-1}{2k}+
\sum_{m\ge 1} \genfrac{(}{)}{0pt}{}{n-m-1}{2k-3}=
\genfrac{(}{)}{0pt}{}{n}{2k-1}+\genfrac{(}{)}{0pt}{}{n-1}{2k+1}+
\genfrac{(}{)}{0pt}{}{n-1}{2k-2}\ .
\]
(iii) \ This follows from (i) and (ii) via a straightforward calculation.
\end{proof}

\begin{thm}  \label{thm:covD}
If\/ $\g=\sone$, i.e., $\Delta$ is of type $\GR{D}{n}$, then
\begin{itemize}
\item[\sf (i)] \ $\cku_\Ab(q)=\sum_{k\ge 0}\bigl(\genfrac{(}{)}{0pt}{}{n+2}{2k}-
4\genfrac{(}{)}{0pt}{}{n-1}{2k-2}\bigr )q^k=\sum_{k\ge 0}\bigl(
\genfrac{(}{)}{0pt}{}{n}{2k}+\genfrac{(}{)}{0pt}{}{n-1}{2k-1}+
\genfrac{(}{)}{0pt}{}{n-2}{2k-1}+\genfrac{(}{)}{0pt}{}{n-2}{2k-4}\bigr)q^k$;
\item[\sf (ii)] \ $\ckl_\Ab(q)=\sum_{k\ge 0}\bigl(\genfrac{(}{)}{0pt}{}{n}{2k+1}{+}
\genfrac{(}{)}{0pt}{}{n}{2k-2}\bigr) q^k$;
\item[\sf (iii)] \ $\gD_\Ab(q)=-\sum_{k\ge 0}
\bigl(\genfrac{(}{)}{0pt}{}{n-2}{2k+1}+\genfrac{(}{)}{0pt}{}{n-3}{2k}\bigr)q^k$.
\end{itemize}
\end{thm}
\begin{proof}
Here the set of commutative roots is 
\[
   \{\esi_i+\esi_j\mid 1\le i<j \le n\} \cup
\{\esi_1-\esi_i\mid 2\le i\le n\} \cup \{\esi_i-\esi_n\mid 2\le i\le n-1\} .
\]
This set is represented by a skew Ferrers diagram
with row lengths $(2n-2,n-1,n-2,\ldots,2)$. See the sample figure for $\Delta$ of type
$\GR{D}{6}$,
where $\esi_{ij}=\esi_i-\esi_j$ and $\ov{\esi}_{ij}=\esi_i+\esi_j$.
%where the leftmost (resp. rightmost) box is $\ap_1=\esi_1-\esi_2$
%(resp. $\theta=\esi_1+\esi_2$). 
%The two boxes in the lowest row are $\ap_{i-1}=\esi_{n-1}-\esi_n$
%and $\ap_{i}=\esi_{n-1}+\esi_n$, respectively.
The roots $\esi_1-\esi_2,\ldots,\esi_1-\esi_{n-1}$ form the {\it tail\/} of the Ferrers diagram for the set of commutative roots
in type $\GR{D}{n}$.

\begin{figure}[h]
\begin{center}
\setlength{\unitlength}{0.02in}
\Yboxdim3.4ex  
{\young(\zab\zac\zad\zae\zaf\zafn\zaen\zadn\zacn\zabn ,::::\zbf\zbfn\zben\zbdn\zbcn  ,::::\zcf\zcfn\zcen\zcdn  ,::::\zdf\zdfn\zden  ,::::\zef\zefn)}
\caption{The (po)set of commutative roots 
for $\mathfrak{so}_{12}$} \label{ab:sone}
\end{center}
\end{figure}

\noindent
To a great extent, convention from the proof  of Theorem~\ref{thm:covB} apply here.
%We follow here the same convention as in the proof of Theorem~\ref{thm:covB}. 
However, the notable distinction of this diagram from Figure~\ref{ab:sono} is that here
the pairs of roots in the two central columns  (i.e., $\esi_{in}, \ov{\esi}_{in}$ for $i<n$)
are incomparable.

\noindent
(i) \ Let us compute the number of abelian ideals $\ah$ with $k$ generators.

$(\heartsuit_1)$ \ Suppose that $\ah$ has the tail of length $m$, $m\ge 1$, i.e.,
$\ah$ has the generator $\esi_1-\esi_{n-m}$ in the upper row.
Then the rest of this row (to the right) is also in the ideal, and
the ideal is determined by its part lying in rows from 2 to $n-1$. 
The condition of being abelian means that
$\ah$ cannot have elements from $m+2$ leftmost columns in this lower part.
(Formally: if $\esi_1-\esi_{n-m}\in \ah$, then for all other roots $\esi_i+\esi_j
\in\ah$, $2\le i\le j$, we must have $j\le n-1-m$.) Hence as a degree of freedom 
for further constructing $\ah$
we obtain a triangle of size $n-3-m$, where an ideal with $k-1$ generators
has to be chosen. The symplectic case shows that the number of possibilities equals
$\genfrac{(}{)}{0pt}{}{n-m-2}{2k-2}$. Thus, in this case we have
$\sum_{m\ge 1}\genfrac{(}{)}{0pt}{}{n-m-2}{2k-2}=\genfrac{(}{)}{0pt}{}{n-2}{2k-1}$
possibilities.

$(\heartsuit_2)$ \ Suppose that an ideal $\ah$ has no tail, i.e., $\esi_1-\esi_{n-1}\not\in \ah$.
Consider all the relevant variants.
\begin{enumerate}
\item $\esi_1-\esi_n, \esi_1+\esi_n\in \ah$. \ These two roots are generators of 
$\ah$, so that we have to choose an ideal with $k-2$ generators in the triangle of
size $n-3$. This yields $\genfrac{(}{)}{0pt}{}{n-2}{2k-4}$ possibilities. 
\item $\esi_1-\esi_n\not\in \ah$. \ Then have to choose an ideal with $k$ 
generators in the triangle of
size $n-1$. This yields $\genfrac{(}{)}{0pt}{}{n}{2k}$ possibilities. 
\item $\esi_1+\esi_n\not\in \ah$. \ This part is the same as previous one, and we
obtain $\genfrac{(}{)}{0pt}{}{n}{2k}$ possibilities.
\item  In items (2) and (3), we have counted twice the ideals that
contain neither $\esi_1-\esi_n$ nor $\esi_1+\esi_n$, i.e., the ideals with $k$ generators
that fit in the triangle of size $n-2$. Therefore $\genfrac{(}{)}{0pt}{}{n-1}{2k}$ must
be subtracted. 
\end{enumerate}
Thus, if $\ah$ has no tail, one obtains the sum
$\genfrac{(}{)}{0pt}{}{n-2}{2k-4}+\genfrac{(}{)}{0pt}{}{n}{2k}+
\genfrac{(}{)}{0pt}{}{n-1}{2k-1}$.
\\
Combining $(\heartsuit_1)$ and $(\heartsuit_2)$ yields the coefficients of $q^k$ presented as 
the sum of four summands. It is a good exercise to transform this sum into the second 
expression in the formulation.
\\[.6ex]
(ii) Counting the ideals with $k$ extensions is even more tedious.
Our approach yields 6 cases and 10 binomial coefficients, which sum luckily up to the two summands in the 
formulation. We only list all the possibilities for the Ferrers diagram
and the corresponding  number of ideals:
\begin{comment}
%%%%%%%%%%%%%%%%%
$\begin{array}{llll}
\bullet & \text{the length of the first row is \ ${\le} n{-}3$}   & - & {\rule{0pt}{2.2ex}}
\genfrac{(}{)}{0pt}{}{n-2}{2k-1} 
\\
\bullet & \text{the length of the first row is \ $n{-}2$}  & - & {\rule{0pt}{2.5ex}}
 \genfrac{(}{)}{0pt}{}{n-2}{2k-4} 
\\
\bullet & \esi_1-\esi_n\not\in \ah$, $\esi_1+\esi_n\in \ah  & - & {\rule{0pt}{2.5ex}}
\genfrac{(}{)}{0pt}{}{n-2}{2k}+\genfrac{(}{)}{0pt}{}{n-2}{2k-3} 
\\
\bullet & \esi_1-\esi_n\in \ah$, $\esi_1+\esi_n\not \in \ah  & - & {\rule{0pt}{2.5ex}}
\genfrac{(}{)}{0pt}{}{n-2}{2k}+\genfrac{(}{)}{0pt}{}{n-2}{2k-3} 
\\
\bullet & \esi_1\pm\esi_n\in \ah   & - & {\rule{0pt}{2.5ex}}
\genfrac{(}{)}{0pt}{}{n-3}{2k}+\genfrac{(}{)}{0pt}{}{n-3}{2k-3}  
\\
\bullet & \text{The ideal has the tail (of length ${\ge} 1$)}  & - & {\rule{0pt}{2.5ex}}
\genfrac{(}{)}{0pt}{}{n-3}{2k+1}+\genfrac{(}{)}{0pt}{}{n-3}{2k-2}. 
\end{array}$
%%%%%%%%%%%%%
\end{comment}
%
\begin{tabbing}
\textbullet\quad \= The ideal has the tail (of length ${\ge} 1$) \quad  \= --\quad \= $\genfrac{(}{)}{0pt}{}{n-3}{2k+1}+\genfrac{(}{)}{0pt}{}{n-3}{2k-2}$. \kill
\textbullet \>  the length of the first row is \ ${\le} n{-}3$   \>   -- \>   
$\genfrac{(}{)}{0pt}{}{n-2}{2k-1}$  {\rule{0pt}{2.5ex}}
\\
\textbullet \>   the length of the first row is \ $n{-}2$  \>   -- \>   {\rule{0pt}{2.5ex}}
$\genfrac{(}{)}{0pt}{}{n-2}{2k-4}$ 
\\
\textbullet \>   $\esi_1-\esi_n\not\in \ah$, $\esi_1+\esi_n\in \ah$  \>   -- \>   {\rule{0pt}{2.5ex}}
$\genfrac{(}{)}{0pt}{}{n-2}{2k}+\genfrac{(}{)}{0pt}{}{n-2}{2k-3}$ 
\\
\textbullet \>   $\esi_1-\esi_n\in \ah$, $\esi_1+\esi_n\not \in \ah$  \>   -- \>   {\rule{0pt}{2.5ex}}
$\genfrac{(}{)}{0pt}{}{n-2}{2k}+\genfrac{(}{)}{0pt}{}{n-2}{2k-3}$ 
\\
\textbullet \>   $\esi_1\pm\esi_n\in \ah$   \>   -- \>   {\rule{0pt}{2.5ex}}
$\genfrac{(}{)}{0pt}{}{n-3}{2k}+\genfrac{(}{)}{0pt}{}{n-3}{2k-3}$ 
\\
\textbullet \>   The ideal has the tail (of length ${\ge} 1$)  \>   -- \>   {\rule{0pt}{2.5ex}}
$\genfrac{(}{)}{0pt}{}{n-3}{2k+1}+\genfrac{(}{)}{0pt}{}{n-3}{2k-2}$. 
\end{tabbing}
(iii) \ This follows from (i) and (ii) via a straightforward calculation.
\end{proof}%

\noindent
Finally, we present the table with complete information about the covering
and deviation polynomials for $\Ab(\g)$.

%%To guess the explicit formula for $\cku$ in the $\GR{D}{n}$-case, we used computer 
%%computations made by G.\,R\"ohrle for $\GR{D}{n}$, $n=6,7,8$.) 
%
\begin{table}[htb]     
\begin{tabular}{l|lll}
\phantom{q}$\g$ & \phantom{quq}$\cku_{\Ab(\g)}$ & 
\phantom{quq}$\ckl_{\Ab(\g)}$ & \phantom{quq} $-\gD_{\Ab(\g)}$ 
\\ \hline
$\GR{A}{n}$ & $\underset{k\ge 0}{\sum}\genfrac{(}{)}{0pt}{}{n+1}{2k}q^k$ & 
$\underset{k\ge 0}{\sum}
\bigl(\genfrac{(}{)}{0pt}{}{n}{2k+1}+\genfrac{(}{)}{0pt}{}{n}{2k-2}\bigr)q^k$
{\rule{0pt}{2.8ex}}  &
$\underset{k\ge 0}{\sum}\genfrac{(}{)}{0pt}{}{n-1}{2k+1}q^k$
\\
$\GR{B}{n}$ &  $\underset{k\ge 0}{\sum}\genfrac{(}{)}{0pt}{}{n+1}{2k}q^k$ & 
$\underset{k\ge 0}{\sum}\bigl(\genfrac{(}{)}{0pt}{}{n-1}{2k+1}{+}
\genfrac{(}{)}{0pt}{}{n}{2k-1}{+}\genfrac{(}{)}{0pt}{}{n-1}{2k-2}\bigr)q^k$
%%$\sum_{k\ge 0}(\genfrac{(}{)}{0pt}{}{n-1}{2k+1}{+}
%%\genfrac{(}{)}{0pt}{}{n+1}{2k-1}{-}\genfrac{(}{)}{0pt}{}{n-1}{2k-3})q^k$
{\rule{0pt}{2.8ex}} 
& $\underset{k\ge 0}{\sum}\genfrac{(}{)}{0pt}{}{n-2}{2k+1}q^k$ 
\\
$\GR{C}{n}$ {\rule{0pt}{2.8ex}}
&  {$\underset{k\ge 0}{\sum}\genfrac{(}{)}{0pt}{}{n+1}{2k}q^k$} &
{$\underset{k\ge 0}{\sum}\genfrac{(}{)}{0pt}{}{n+1}{2k}q^k$} & $0$ 
\\
$\GR{D}{n}$ & 
$\underset{k\ge 0}{\sum}\bigl(\genfrac{(}{)}{0pt}{}{n+2}{2k}-
4\genfrac{(}{)}{0pt}{}{n-1}{2k-2}\bigr )q^k$ {\rule{0pt}{3ex}} &
$\underset{k\ge 0}{\sum}\bigl(\genfrac{(}{)}{0pt}{}{n}{2k+1}+
\genfrac{(}{)}{0pt}{}{n}{2k-2}\bigr)q^k$ & 
$\underset{k\ge 0}{\sum}\bigl(\genfrac{(}{)}{0pt}{}{n-2}{2k+1}{+}
\genfrac{(}{)}{0pt}{}{n-3}{2k}\bigr)q^k$ 
\\
$\GR{E}{6}$ & $1{+}25q{+}\phantom{1} 27q^2{+}11q^3$  {\rule{0pt}{2.4ex}} &
$6{+}21q{+}20q^2{+}17q^3$ & $5{+}\phantom{1} 6q$ 
\\
$\GR{E}{7}$ & $1{+}34q{+}\phantom{1} 60q^2{+}30q^3{+}\phantom{1} 3q^4$ 
& $7{+}35q{+}40q^2{+}43q^3{+}\phantom{1} 3q^4$ 
{\rule{0pt}{2.3ex}} & $6{+}13q$ 
\\
$\GR{E}{8}$ & $1{+}44q{+}118q^2{+}76q^3{+}17q^4$ {\rule{0pt}{2.3ex}} &
$8{+}49q{+}87q^2{+}95q^3{+}17q^4$ & $7{+}19q$  
\\
$\GR{F}{4}$ & $1{+}10q{+}\phantom{11} 5q^2$ & 
$2{+}\phantom{1} 8q{+}\phantom{1}6q^2$ & $1$ 
\\
$\GR{G}{2}$ & %\multicolumn{2}{c}%
{$1{+}\phantom{1} 3q$}  & {$1{+}\phantom{1} 3q$} & $0$ 
\\ \hline
\end{tabular} 
\vskip.7ex 
\caption{ The covering and deviation
polynomials for $\Ab(\g)$}  \label{table:ab}
\end{table}%
\vskip.5ex\noindent
{\it Some observations related to Table~\ref{table:ab}}  \nopagebreak

1. For $\GR{A}{2n}$, we have
$\deg\ckl_\Ab-\deg\cku_\Ab=1$. For all other cases the degrees are equal.
Perhaps, the reason is that the Coxeter number is odd if and only if $\Delta$ is of type
$\GR{A}{2n}$.

2. One may observe that there are several regularities
in Table~\ref{table:ab}.
For all classical series, both  covering polynomials 
satisfy the recurrence relation 
\begin{equation}  \label{recurs}
 \ck_{\Ab(\GR{X}{n})}(q)=2\ck_{\Ab(\GR{X}{n-1})}(q)+(q-1)\ck_{\Ab(\GR{X}{n-2})}(q)\ , 
\end{equation}
where $\GR{X}{}\in \{\mbox{\bf A,B,C,D}\}$ and
$\ck$ is either $\cku$ or $\ckl$. 
Furthermore, the sequence $\GR{E}{3}=\GR{A}{2}\times \GR{A}{1}$, 
$\GR{E}{4}=\GR{A}{4},\GR{E}{5}=\GR{D}{5}$, $\GR{E}{6}$,
$\GR{E}{7},\GR{E}{8}$
can be regarded as the `exceptional' series, and for this series the same
recurrence relation is satisfied for $\cku$ (but not for $\ckl$). 
Here one might wonder about the possible meaning of the  polynomial $\cku$
corresponding to ``$\GR{E}{9}$''.

It also follows from \eqref{recurs}
that $\ck_{\Ab(\GR{X}{n})}(1)=2\ck_{\Ab(\GR{X}{n-1})}(1)$, which "explains" the equality
$\#\Ab(\GR{X}{n})=2^n$.

3. The upper covering polynomials for $\GR{A}{n}, \GR{B}{n}, \GR{C}{n}$ are the 
same. But the lower covering polynomial distinguishes these series.
Furthermore, if the Dynkin diagram has no branching nodes, then $\cku_\Ab(\g)$
depends only on $n$. That is, the upper covering polynomial for $\GR{F}{4}$ 
(resp. $\GR{G}{2}$) is equal
to that for $\GR{A}{4}$ (resp. $\GR{A}{2}$).

On the other hand, the lower covering polynomials are equal for
$\GR{A}{n}$ and $\GR{D}{n}$, and the deviation polynomial for $\GR{B}{n}$ is equal 
to that for $\GR{A}{n-1}$.
\\[1ex]
It would be interesting to find an explanation of these coincidences.
%%(a bijective proof?).

{\bf Remark.} For the non-reduced root system $\GR{BC}{n}$, one can also consider 
combinatorial abelian ideals. However, these are exactly the same as in the 
symplectic case.

%%%%%%%%%%%%%%%%   Section 6

\section{Some questions and open problems}
\label{some_spec}

\noindent
$1^o$. Is there a combinatorial interpretation of the values $\cku_\PP(-1)$ and
$\ckl_\PP(-1)$? Specifically, for the posets $\AD$, $\AD_0$, and $\Ab$, one might expect
some characteristics of the corresponding root system.
%(At least, if $\PP$ is a modular or distributive lattice.)
\\[.6ex]
$2^o$. Various examples considered in the paper show that in many cases
the deviation polynomial of a poset has the nonzero coefficients of the same sign.
It would be interesting to find a general pattern for this phenomenon.
Of course, it is not always the case. For instance,
if $\gD_{\PP_1}$ (resp. $\gD_{\PP_2}$) has positive (resp. negative) coefficients,
then Lemma~\ref{prod}(ii) shows that the deviation
polynomial of $\PP_1\times \PP_2$ may have coefficients of both signs.
It is not hard to produce a concrete example.
However, I conjecture that the following is true: 

%%{\it Suppose $\PP$ is a graded meet-semilattice. Then
%%the nonzero coefficients of $\gD_\PP$ are negative.}
%%\\[.7ex]
%%At least, we can prove that in this case $\gD_\PP(1)\le 0$.
%%It is also worth stating a particular case of this conjecture:

{\it Suppose that $\PP=J^*(\LL)$, and let $\PP({\le} m)$ be the subposet of upper ideals 
of $\LL$ whose cardinality is at most $m$. Then 
$\gD_{\PP({\le} m)}$ has non-positive coefficients for all $m$.}
\\[.6ex]
In the special case of $J^*(\LL)\setminus \{\hat 1\}$, this is verified in 
Example~\ref{ex:transform}.
\\[.6ex]
\noindent
$3^o$. Using Table~\ref{table:ab}, one can compute the values $-\gD_\Ab(1)$.
In the serial cases, these values are quite simple: 

$2^{n-2}$ for $\GR{A}{n}$, \ $n\ge 2$; $2^{n-3}$ for $\GR{B}{n}$, $n\ge 3$; \ 
$0$ for $\GR{C}{n}$, $n\ge 1$; \ 
$2^{n-3}+2^{n-4}$ for $\GR{D}{n}$, $n\ge 4$.
\\
Hopefully, there could be a uniform general description for them.
One may notice that if the Dynkin diagram has no branching nodes, then
this value equals $2^m$ with $m=\#(\Pi_l)-2$; or $0$, if
$\#(\Pi_l)=1$. (The case of $\GR{G}{2}$ and $\GR{F}{4}$ is included here.)
But I have no idea how to explain the values for series $\bf D$ and $\bf E$.


\begin{thebibliography}{Pa95}

\bibitem{ath1} {\sc C.A.~Athanasiadis}. On a refinement of the generalized 
Catalan numbers for Weyl groups, {\it Trans. AMS}, {\bf 357}(2005), 179--196.

%%Generalized Catalan numbers, 
%%Weyl groups and arrangements of hyperplanes, {\it Bull. London Math. Soc.},
%%{\bf 36}(2004), 294--302.

\bibitem{bessis} {\sc D.~Bessis}. The dual braid monoid,
{\it Ann. Sci. \'Ecole Norm. Sup.}, S\'er. IV, {\bf 36}\,(2003), 647--683. 

\bibitem{bour} {\sc N.~Bourbaki.}
"Groupes et alg\`ebres de Lie", Chapitres 4,5 et 6,
Paris: Hermann 1975.

\bibitem{cp1} {\sc P.~Cellini} and {\sc P.~Papi}.
ad-nilpotent ideals of a Borel subalgebra, {\it J. Algebra},  {\bf 225}\,(2000),  130--141.

\bibitem{cp2} {\sc P.~Cellini} and {\sc P.~Papi}.
ad-nilpotent ideals of a Borel subalgebra II, {\it J. Algebra}, {\bf 258}\,(2002), 112--121.

%\bibitem{cp3} {\sc P.~Cellini} and {\sc P.~Papi}.
%Abelian ideals of Borel subalgebras and affine Weyl groups, 
%{\it Adv. Math.}, {\bf 187}(2004), 320--361.

\bibitem{chap02} {\sc F.~Chapoton}. Data file on cluster graphs, 
{\it Manuscript dated} August 30, 2002.

\bibitem{chap04} {\sc F.~Chapoton}. Enumerative properties of generalized associahedra. 
{\it S\'em. Lothar. Combin.} {\bf 51}(2004), Art. B51b, 16 pp. (electronic).

%\bibitem{djok-et}
%{\sc D.~Djokovi\'c,  P.~Check} and {\sc  J.-Y.~H\'ee}.
%On closed subsets of root systems, 
%{\it Canad. Math. Bull.} {\bf 37}(1994), no. 3, 338--345.  

\bibitem{dilw}
{\sc R.P.~Dilworth}. Proof of a conjecture on finite modular lattices,
{\it Annals Math.}, {\bf 60}\,(1954),  359--364.

%%\bibitem{cluster} {\sc S.V.~Fomin} and {\sc A.V.~Zelevinsky}.
%%$Y$-system and generalized associahedra, 
%%{\it Ann. of Math.}, {\bf 158}(2003), no. 3, 977--1018.

\bibitem{hump}
{\sc J.E.~Humphreys}. ``Reflection Groups and Coxeter Groups'',
Cambridge Univ. Press, 1992.

\bibitem{jos}  {\sc A.~Joseph}. Orbital varieties of the minimal orbit
{\it Ann. Scient. \'Ec. Norm. Sup. (4)},  {\bf 31}\,(1998), 17--45.

\bibitem{ko1} {\sc B.~Kostant}. The set of abelian ideals of a Borel
subalgebra, Cartan
decompositions, and discrete series representations,
{\it Intern. Math. Res. Notices},  (1998),  no.\,5, 225--252.

\bibitem{kratt} {\sc C.~Krattenthaler}. The $F$-triangle of the generalised cluster complex,
in: ``Topics in discrete mathematics'', 93--126, Algorithms Combin., 26, Springer, Berlin, 2006.

\bibitem{aif99} {\sc D.~Panyushev}.
On spherical nilpotent orbits and beyond, 
{\it Ann. Inst. Fourier}, {\bf 49}\,(1999), 1453--1476.

\bibitem{imrn} {\sc D.~Panyushev}. Abelian ideals of a Borel subalgebra and
long positive roots, {\it Intern. Math. Res. Notices\/} (2003), no.\,35,
1889--1913.

\bibitem{duality} {\sc D.~Panyushev}. {\sf ad}-nilpotent ideals of a
Borel subalgebra: generators and duality, {\it J. Algebra}, 
{\bf 274}\,(2004), 822--846.

%%\bibitem{long} {\sc D.~Panyushev}. Long abelian ideals, {\it Adv. Math.},
%%{\bf 186}(2004), 307--316.

\bibitem{rodstv} {\sc D.~Panyushev}. The poset of positive roots and its
relatives, {\it J. Algebraic Combin.}, {\bf 23}\,(2006), 79--101.

\bibitem{mmj} {\sc D.~Panyushev}. Properties of weight posets for weight multiplicity 
free representations, {\it Moscow Math. J.}, {\bf 9}, no.\,4 (2009),  867--883.

\bibitem{pr}  {\sc D.~Panyushev} and {\sc G.~R\"ohrle}.
Spherical orbits and abelian ideals, {\it Adv. Math.}, {\bf 159}\,(2001),
229--246.

\bibitem{charney}  {\sc V.~Reiner} and {\sc V.~Welker}.
On the Charney-Davis and Neggers-Stanley conjectures, 
{\it J. Combin. Theory. Ser.~A} \ {\bf 109}\,(2005), no.\,2, 247--280.

\bibitem{R} {\sc G.~R\"ohrle}. 
On normal abelian subgroups of parabolic groups,
{\it Ann. Inst. Fourier}, {\bf 48}\,(1998), 1455--1482.

\bibitem{eric} {\sc E.~Sommers.} $B$-stable ideals in the nilradical of a
Borel subalgebra, {\it Canad. Math. Bull.}, {\bf 48}\,(2005),  460--472.

\bibitem{rstan1} {\sc R.P.~Stanley}.
``Enumerative Combinatorics'',  vol. 1. Cambridge Univ. Press, 1997.

\bibitem{stembr} {\sc J.R.~Stembridge}. Trapezoidal chains and antichains,
{\it Europ. J. Combin.}, {\bf 7}\,(1986), 377--387.


\end{thebibliography}
\end{document}